\documentclass[12pt]{article}
\usepackage{amsfonts}

\usepackage[all]{xy}

\def\squarebox#1{\hbox to #1{\hfill\vbox to #1{\vfill}}}
\newcommand{\qed}{\hspace*{\fill}
\vbox{\hrule\hbox{\vrule\squarebox{.667em}\vrule}\hrule}\smallskip}
\newtheorem{teorema}{Theorem}[section]
\newtheorem{lema}[teorema]{Lemma}
\newtheorem{corolario}[teorema]{Corollary}
\newtheorem{proposicao}[teorema]{Proposition}
\newtheorem{defi}[teorema]{Definition}
\newenvironment{profe}{\noindent {\bf Proof:}}{\hfill $\qed $ \newline}
\newenvironment{remarque}{\noindent {\bf Remark:}}{\vspace{0.5cm}}
\newenvironment{exemplo}{\noindent {\bf Example:}}{}
\begin{document}

\title{Conley index and stable sets \\
for flows on flag bundles}

\author{Mauro Patr\~{a}o\footnote{Departamento de Matem\'{a}tica, 
Universidade de Bras\'{\i}lia. Campus Darcy Ribeiro. Bras\'{\i}lia, DF, Brasil.
 \textit{e-mail: mpatrao@mat.unb.br }.} 
\and Luiz A. B. San Martin\footnote{Departamento de Matem\'{a}tica. 
Universidade Estadual de Campinas. Cx. Postal 6065. 13.083-859 Campinas, SP,
Brasil. \textit{e-mail: smartin@ime.unicamp.br}
Supported by CNPq grant n$^{\circ }$ $305513/2003$-$6$ and FAPESP grant 
n$^{\circ }$ $02/10246$-$2$. } \\
Campinas SP, Brasil \and Lucas Seco\footnote{Departamento de Matem\'{a}tica. 
Universidade Estadual de Campinas. Cx. Postal 6065. 13.083-859 Campinas, SP,
Brasil. \textit{e-mail: ajsantana@uem.br}.
Supported by FAPESP grant n$^{\protect\underline{\circ }}$ 07/52390-6 } }

\maketitle

\begin{abstract}
Consider a continuous flow of automorphisms of a $G$-principal bundle which
is chain transitive on its compact Hausdorff base. Here $G$ is a connected
noncompact semi-simple Lie group with finite center. The finest Morse
decomposition of the induced flows on the associated flag bundles were
obtained in previous articles. Here we describe the stable sets of these
Morse components and, under an additional assumption, their Conley indices.
\end{abstract}



\section{Introduction}

The subject matter of this paper are the flows on flag bundles. The
questions to be treated are the Conley indices and the stable (and unstable)
sets of the Morse components.

We consider the following setup (that already appears in \cite{smbflow},
\cite{msm}, \cite{ll1} and elsewhere): Start with a principal bundle $%
Q\rightarrow X$ whose structural group $G$ is semi-simple non compact (or
slightly more generally the Lie algebra of $G$ is reductive and $G$
satisfies some mild assumptions) and $X$ is a compact Hausdorff space. The
group $G$ acts freely (on the right) on $Q$ having $X$ as the space of
orbits. We assume throughout that the bundle is locally trivial. Now, let $%
\phi _{t}$ ($t\in \Bbb{Z}$ or $\Bbb{R}$) be a flow on $Q$ which commutes
with the right action of $G$: $\phi _{t}\left( p\cdot a\right) =\phi
_{t}\left( p\right) \cdot a$, $p\in Q$ and $a\in G$, where $p\cdot a$
denotes the right action. Then $\phi _{t}$ induces flows on any associated
bundle $Q\times _{G}F$ where $F$ is a space acted on the left by $G$. We
mainly take $F$ to be a flag manifold of $G$.

This set up includes the linear flows on vector bundles and on the
corresponding projective and Grassmann bundles, for which there is an
extensive literature (we cite, for instance, Colonius-Kliemann \cite{ck},
Sacker-Sell \cite{sacsell}, Salamon-Zehnder \cite{sz}, and their
references). For these bundles one can take $G=\mathrm{Gl}\left( n,\Bbb{R}%
\right) $ for the structure group.

The finest Morse decomposition for the flows on flag bundles were described
in \cite{smbflow} and \cite{msm} with the assumption that the flow on the
base space is chain transitive. From this description it emerges a reduction
of $Q$ to a subbundle $Q_{\phi }$ that is invariant under the flow. The
structure group of the reduced bundle is the centralizer $Z_{H_{\phi
}}\subset G$ of a suitable element $H_{\phi }$ of the Lie algebra $\frak{g}$
of $G$. When $G=\mathrm{Gl}\left( n,\Bbb{R}\right) $ this centralizer is a
subgroup of block diagonal matrices, which suggested the name of block
reduction of the flow (see \cite{ll1}, Section 5, and Section \ref{secstable}
below).

The block reduction is one of our main tools. We use it to transport to flag
bundles constructions made on the flag manifolds. In fact, by the very
definition of $Q_{\phi }$ in \cite{ll1} the union of Morse components in a
flag bundle is the set $\{q\cdot \mathrm{fix}_{H_{\phi }}:q\in Q_{\phi }\}$
where $\mathrm{fix}_{H_{\phi }}$ is the set of fixed points of the action on
a flag manifold of $\exp H_{\phi }$ (below we use a different notation, that
distinguish the Morse components as well as the specific flag manifold). In
this paper we make a similar construction for the stable and unstable sets
of the Morse components. They are obtained by plugging into the fibers of a
flag bundle the stable sets on the flag manifolds (see Theorem \ref
{teo:var-estavel}).

The stable and unstable sets are used to describe the homotopy (and hence
the cohomology) Conley index of a Morse component $\mathcal{M}$. It turns
out that the homotopy index is the Thom space of a vector bundle $\mathcal{U}
$ over $\mathcal{M}$ (see Theorem \ref{teo:indice-conley}; we recall de
definition of the Thom space in Section \ref{preliminar}). The fibers of $%
\mathcal{U}$ are the tangent spaces (along the the fibers of the flag
bundle) to the unstable set of $\mathcal{M}$. By Thom isomorphism the
(co)homology Conley index is the displaced (co)homology of $\mathcal{M}$. In
some cases (e.g. the base space is contractible) the (co)homology of $%
\mathcal{M}$ reduces to the (co)homology of a flag manifold (see \cite
{kocherlakota} and \cite{gelfand}).

To prove Theorem \ref{teo:indice-conley} we first linearize the flow around $%
\mathcal{M}$ by building a conjugation between a linear flow on $\mathcal{U}
$ and the flow restricted to a neighborhood of $\mathcal{M}$.

However, in order that our construction of the linearization become well
defined we make an additional assumption about the flow. Namely, that the
block reduction subbundle $Q_{\phi }$ admits a further $\phi _{t}$-invariant
reduction to a subbundle $C_{\phi }$ whose structural group $L_{H_{\phi }}$
is a subgroup of $Z_{H_{\phi }}$. The subgroup $L_{H_{\phi }}$ is the direct
product of a compact group by an abelian vector group. If $G=\mathrm{Gl}%
\left( n,\Bbb{R}\right) $ then $L_{H_{\phi }}$ a subgroup of block diagonal
matrices where each block is a conformal matrix (product of an orthogonal
matrix by a scalar one). This suggests to call $C_{\phi }$ a block conformal
reduction of the flow.

In general such a conformal reduction may not exist. Although we do not know
the full significance of this geometric assumption there are some clues in
the literature about its reliability. At this regard we mention the results
on Jordan decomposition of cocycles by Arnold-Cong-Oseledets \cite{aco},
where block conformal matrices are produced by the support of invariant
measures on projective and sphere bundles. On the other hand our $Z_{H_{\phi
}}$-reduction $Q_{\phi }$ is produced by the chain recurrent components on
flag bundles. Hence one may suspect that the reason for the existence of a
conformal reduction $C_{\phi }$ stays in the relation between the supports
of invariant measures and the chain components.

\section{Preliminaries\label{preliminar}}

Troughout the article we will use the notation of principal and associated
bundles of Kobayashi-Nomizu \cite{kn}. In what follows we establish some
other notations and preliminar results.

Let $\pi :Q\to X$ be a principal bundle with structural group $G$ and
compact Hausdorff base space $X$. Let $\phi _{t}$ be a right invariant flow
on $Q$. We assume throughout that the induced flow on $X$ is chain
transitive.

\subsection{Flows on topological spaces\label{preliminar-conley}}

Let $\phi :\mathbb{T}\times E\to E$ be a continuous flow on a compact
Hausdorff space $E$, with discrete $\mathbb{T}={\mathbb
Z}$ or continuous $\mathbb{T}={\mathbb R}$ time. Fix an invariant set $%
\mathcal{M}\subset E$. We define its stable and unstable sets respectively
as
\[
\mathrm{st}(\mathcal{M})=\{x\in E:\omega (x)\subset \mathcal{M}\},\qquad
\mathrm{un}(\mathcal{M})=\{x\in E:\omega ^{*}(x)\subset \mathcal{M}\}.
\]
and its attractor and repeller domains respectively as
\[
\mathcal{A}\left( \mathcal{M}\right) =\{x\in E:\Omega (x)\subset \mathcal{M}%
\},\qquad \mathcal{R}\left( \mathcal{M}\right) =\{x\in E:\Omega^*(x)\subset
\mathcal{M}\},
\]
where $\omega(x)$, $\omega^*(x)$ are the limit sets of $x$, and $\Omega(x)$,
$\Omega^*(x)$ are the chain limit sets of $x$. Since we have that $\omega(x)
\subset \Omega(x)$, it follows that $\mathrm{st}(\mathcal{M}) \subset
\mathcal{A}\left( \mathcal{M}\right)$ and analogously $\mathrm{un}(\mathcal{M%
}) \subset \mathcal{R}\left( \mathcal{M}\right)$. For the concepts of Morse
decompositions and its relations with chain transitivity see \cite{c}, \cite
{ck}, \cite{patr}. If $\{\mathcal{M}_{i}\}_{i\in I}$ is a Morse
decomposition of $E$ then $E$ decomposes as the disjoint union of stable
sets $\mathrm{st}(\mathcal{M}_{i})$.

We now recall the definition of the Conley index. We note that, when working
with the Conley index, we restrict ourselves to the case of a
continuous-time flow. A neighborhood $U \subset E$ of the invariant set $%
\mathcal{M}$ is isolating if $\mathcal{M}$ is the maximal invariant set
inside $U$. A pair $(N_1, N_0)$ of subsets of $E$ with $N_1 \supset N_0$ is
an index pair for $\mathcal{M}$ if it satisfies the following three
conditions

\begin{itemize}
\item[1)]  $\mathrm{cl}(N_{1}\backslash N_{0})$ is an isolating neighborhood
of $\mathcal{M}$.

\item[2)]  $N_{0}$ is positively invariant with respect to $N_{1}$, i.e., if
$x\in N_{0}$ and $\phi _{t}(x)\in N_{1}$, for some $t\geq 0$, then $\phi
_{t}(x)\in N_{0}$.

\item[3)]  If $x\in N_{1}$ is such that there exists $t^{\prime }>0$ with $%
\phi _{t^{\prime }}(x)\notin N_{1}$, then there exists $t_{m}\geq 0$ such
that $\phi _{t_{m}}(x)\in N_{0}$ and $\phi _{t}(x)\in N_{1}$, for every $%
t\in [0,t_{m}]$.
\end{itemize}

The Conley index of $\mathcal{M}$ is defined as the pointed homotopy class $%
h(\mathcal{M})$ of the quotient space $N_{1}/N_{0}$ with basepoint $[N_{0}]$%
. A fundamental result of the abstract theory (see Conley-Zehnder \cite{cz})
is that index pairs for an isolated invariant set $\mathcal{M}$ always exist
and that the Conley index of $\mathcal{M}$ is independent of a chosen index
pair.

Let $H^{*}$ denote $\check{\mathrm{C}}$ech cohomology with coefficients in
some fixed ring. The cohomological Conley index of $\mathcal{M}$ is the
reduced cohomology of the Conley index with respect to its basepoint and is
denoted by
\[
CH^{*}(\mathcal{M})=\widetilde{H}^{*}(h(\mathcal{M})).
\]
We now recall a construction due to Thom that will be used in the
computation of our indexes. Let $p:\mathcal{V}\to B$ be an $n$-dimensional
vector bundle over a paracompact Hausdorff base $B$. The Thom space $T(%
\mathcal{V})$ of $\mathcal{V}$ is the pointed homotopy class of the quotient
space $\mathcal{V}/(\mathcal{V}-0)$ with basepoint $[(\mathcal{V}-0)]$. A
Thom class $U$ for $\mathcal{V}$ is an element $U\in H^{n}(\mathcal{V},%
\mathcal{V}-0)$ such that its restriction to $H^{n}(\mathcal{V}_{b},(%
\mathcal{V}-0)_{b})$ is a generator, for all $b\in B$. A Thom class for ${%
\mathbb Z}_{2}$ coefficients always exists, and it exists for ${\mathbb Z}$
coefficients if and only if $\mathcal{V}$ is orientable (see e.g. \cite
{hatcher-at}, Theorem 4D.10). If there exists a Thom class $U$ then the map
\[
H^{*}(B)\to H^{*+n}(\mathcal{V},\mathcal{V}-0),\qquad \alpha \mapsto
p^{*}(\alpha )\cup U.
\]
is Thom isomorphism (see \cite{hatcher-at}, Corollary 4D.9). Here $\cup $ is
the cup product and $p^{*}$ is the induced map in cohomology of the
projection $p$. Passing from relative to reduced cohomology the Thom
isomorphism becomes
\[
H^{*}(B)\to \widetilde{H}^{*+n}(T(\mathcal{V})),\qquad \alpha \mapsto
p^{*}(\alpha )\cup U.
\]
Put a riemannian metric $|\cdot |$ in $\mathcal{V}$, let $D(\mathcal{V}%
)=\{v\in \mathcal{V}:\,|v|\leq 1\}$ be the disk bundle and let $S(\mathcal{V}%
)=\{v\in \mathcal{V}:\,|v|=1\}$ be the sphere bundle of $\mathcal{V}$. The
inclusion of pairs $(D\mathcal{V},S\mathcal{V})\subset (\mathcal{V},\mathcal{%
V}-0)$ is easily seen to induce a pointed homotopy equivalence between $D(%
\mathcal{V})/S(\mathcal{V})$ with basepoint $[S(\mathcal{V})]$ and the Thom
space $T(\mathcal{V})$, so the Thom space of $\mathcal{V}$ can alternatively
be defined as the pointed homotopy type this quotient.

Now we recall the Morse equation obtained in \cite{cz}. Let $(Y_{1},Y_{0})$
be a pair of compact Hausdorff spaces with $Y_{1}\supset Y_{0}$ and denote
by $H(Y_{1},Y_{0})$ its $\check{\mathrm{C}}$ech cohomology with coefficients
in some fixed ring. Assuming each module $H^{j}(Y_{1},Y_{0})$ to be of
finite rank we define the Poincar\'{e} polynomial of the pair $(Y_{1},Y_{0})$
by the formal sum
\[
P(t,Y_{1},Y_{0}):=\sum_{j\geq 0}t^{j}\,\mathrm{rank}\,H^{j}(Y_{1},Y_{0}).
\]
We define the Poincar\'{e} polynomial of the space $Y$ by $%
P(t,Y):=P(t,Y,\emptyset )$ and the Poincar\'{e} polynomial of the index of
an isolated invariant set $\mathcal{M}\subset E$ by
\[
CP(t,\mathcal{M})=P(t,N_{1},N_{0}),
\]
where $(N_{1},N_{0})$ is an index pair for $\mathcal{M}$. Since for $\check{%
\mathrm{C}}$ech cohomology we have $H(N_{1},N_{0})=\widetilde{H}%
(N_{1}/N_{0})=CH(\mathcal{M})$, it follows that
\[
CP(t,\mathcal{M}):=\sum_{j\geq 0}t^{j}\,\mathrm{rank}\,CH^{j}(\mathcal{M}),
\]
so that this polynomial does not depend on the chosen index pair. The Morse
equation relates the cohomology of the indices with the cohomology of the
whole space $E$.

\begin{teorema}
\label{teo:morse-eq-abstrata}Let $\{\mathcal{M}_{i}\}_{i\in I}$ be a Morse
decomposition of $E$, then
\[
\sum_{i\in I}CP(t,\mathcal{M}_{i})=P(t,E)+(1+t)R(t).
\]
If the coefficient ring is $\Bbb{Z}$ then the coefficients of $R(t)$ are
non-negative.
\end{teorema}

\subsection{Semi-simple Lie Theory\label{preliminar-lie}}

For the theory of semi-simple Lie groups and their flag manifolds we refer
to Duistermat-Kolk-Varadarajan \cite{dkv}, Helgason \cite{Helgason}, Knapp
\cite{kpp} and Warner \cite{w}. To set notation let $G$ be a connected
noncompact semi-simple Lie group with Lie algebra $\frak{g}$. We assume
throughout that $G$ has finite center. Fix a Cartan involution $\theta $ of $%
\frak{g}$ with Cartan decomposition $\frak{g}=\frak{k}\oplus \frak{s}$. The
form $B_{\theta }\left( X,Y\right) =-\langle X,\theta Y\rangle $, where $%
\langle \cdot ,\cdot \rangle $ is the Cartan-Killing form of $\frak{g}$, is
an inner product.

Fix a maximal abelian subspace $\frak{a}\subset \frak{s}$ and a Weyl chamber
$\frak{a}^{+}\subset \frak{a}$. We let $\Pi $ be the set of roots of $\frak{a%
}$, $\Pi ^{+}$ the positive roots corresponding to $\frak{a}^{+}$, $\Sigma $
the set of simple roots in $\Pi ^{+}$ and $\Pi^- = - \Pi^+$ the negative
roots. For a root $\alpha \in \Pi$ we denote by $H_\alpha \in \frak{a}$ its
coroot so that $B_\theta(H_\alpha,H) = \alpha(H)$ for all $H \in \frak{a}$.
The Iwasawa decomposition of the Lie algebra $\frak{g}$ reads $\frak{g}=%
\frak{k}\oplus \frak{a}\oplus \frak{n}^{\pm}$ with $\frak{n}%
^{\pm}=\sum_{\alpha \in \Pi ^{\pm}}\frak{g}_{\alpha }$ where $\frak{g}%
_{\alpha }$ is the root space associated to $\alpha $. As to the global
decompositions of the group we write $G=KS$ and $G=KAN^\pm$ with $K=\exp
\frak{k}$, $S=\exp \frak{s}$, $A=\exp \frak{a}$ and $N^{\pm}=\exp \frak{n}%
^{\pm}$.

The Weyl group $\mathcal{W}$ associated to $\frak{a}$ is the finite group
generated by the reflections over the root hyperplanes $\alpha =0$ in $\frak{%
a}$, $\alpha \in \Pi $. $\mathcal{W}$ acts on $\frak{a}$ by isometries and
can be alternatively be given as $\mathcal{W}=M^{*}/M$ where $M^{*}$ and $M$
are the normalizer and the centralizer of $A$ in $K$, respectively. We write
$\frak{m}$ for the Lie algebra of $M$. There is an unique element $w^{-}\in
\mathcal{W}$ which takes the simple roots $\Sigma $ to $-\Sigma $, $w^{-}$
is called the principal involution of $\mathcal{W}$. The Bruhat-Chevalley
order in $\mathcal{W}$ is a partial order given as follows. Take for $w\in
\mathcal{W}$ a reduced expression $w=s_{1}\cdots s_{n}$ as a product of
reflections with respect to the simple roots $\Sigma $. Then $\overline{w}%
\leq w$ if and only if there are integers $1\leq i_{1}\leq \cdots \leq
i_{j}\leq n$ such that $\overline{w}=s_{i_{1}}\cdots s_{i_{j}}$ is a reduced
expression for $\overline{w}$. We have that $w^{-}$ is the greatest element
of $\mathcal{W}$ in this order.

Associated to a subset of simple roots $\Theta \subset \Sigma $ there are
several Lie algebras and groups (cf. \cite{w}, Section 1.2.4): We write $%
\frak{g}\left( \Theta \right) $ for the (semi-simple) Lie subalgebra
generated by $\frak{g}_{\alpha }$, $\alpha \in \Theta $, and put $\frak{k}%
(\Theta )=\frak{g}(\Theta )\cap \frak{k}$, $\frak{a}\left( \Theta \right) =%
\frak{g}\left( \Theta \right) \cap \frak{a}$, and $\frak{n}^{\pm }\left(
\Theta \right) =\frak{g}\left( \Theta \right) \cap \frak{n}^{\pm }$. The
simple roots of $\frak{g}(\Theta )$ are given by $\Theta $, more precisely,
by restricting the functionals of $\Theta $ to $\frak{a}(\Theta )$. The
coroots $H_{\alpha }$, $\alpha \in \Theta $, form a basis for $\frak{a}%
\left( \Theta \right) $. Let $G\left( \Theta \right) $ and $K(\Theta )$ be
the connected groups with Lie algebra $\frak{g}\left( \Theta \right) $ and $%
\frak{k}\left( \Theta \right) $, respectively. Then $G(\Theta )$ is a
connected semi-simple Lie group with finite center. Let $A\left( \Theta
\right) =\exp \frak{a}\left( \Theta \right) $, $N^{\pm }\left( \Theta
\right) =\exp \frak{n}^{\pm }\left( \Theta \right) $. We have the Iwasawa
decomposition $G(\Theta )=K(\Theta )A(\Theta )N^{\pm }(\Theta )$. Let $\frak{%
a}_{\Theta }=\{H\in \frak{a}:\alpha (H)=0,\,\alpha \in \Theta \}$ be the
orthocomplement of $\frak{a}(\Theta )$ in $\frak{a}$ with respect to the $%
B_{\theta }$-inner product and put $A_{\Theta }=\exp \frak{a}_{\Theta }$.
The subset $\Theta $ singles out the subgroup $\mathcal{W}_{\Theta }$ of the
Weyl group which acts trivially $\frak{a}_{\Theta }$. Alternatively $%
\mathcal{W}_{\Theta }$ can be given as the subgroup generated by the
reflections with respect to the roots $\alpha \in \Theta $. The restriction
of $w\in \mathcal{W}_{\Theta }$ to $\frak{a}(\Theta )$ furnishes an
isomorphism between $\mathcal{W}_{\Theta }$ and the Weyl groyp $\mathcal{W}%
(\Theta )$ of $G(\Theta )$

The standard parabolic subalgebra of type $\Theta \subset \Sigma$ with
respect to chamber $\frak{a}^+$ is defined by
\[
\frak{p}_{\Theta }=\frak{n}^{-}\left( \Theta \right) \oplus \frak{m}\oplus
\frak{a}\oplus \frak{n}^{+}.
\]
The corresponding standard parabolic subgroup $P_{\Theta }$ is the
normalizer of $\frak{p}_{\Theta }$ in $G$. It has the Iwasawa decomposition $%
P_{\Theta } = K_\Theta A N^+$. The empty set $\Theta =\emptyset $ gives the
minimal parabolic subalgebra $\frak{p} = \frak{m}\oplus \frak{a}\oplus \frak{%
n}^{+}$ whose minimal parabolic subgroup $P=P_\emptyset $ has Iwasawa
decomposition $P=MAN^+$. Let $\Delta \subset \Theta \subset \Sigma$. Then $%
P_\Delta \subset P_\Theta$, also we denote by $P(\Theta)_\Delta$ the
parabolic subgroup of $G(\Theta)$ of type $\Delta$.

We let $Z_{\Theta }$ be the centralizer of $\frak{a}_{\Theta }$ in $G$ and $%
K_{\Theta }=Z_{\Theta }\cap K$. We have that $K_{\Theta }$ decomposes as $%
K_{\Theta }=MK(\Theta )$ and that $Z_{\Theta }$ decomposes as $Z_{\Theta
}=MG(\Theta )A_{\Theta }$ which implies that $Z_{\Theta }=K_{\Theta
}AN(\Theta )$ is an Iwasawa decomposition of $Z_{\Theta }$ (which is a
reductive Lie group). Let $\Delta \subset \Sigma $, then\footnote{%
Using that $\frak{a}_\Theta = \frak{a}(\Theta)^\perp$ and that $\frak{a}%
(\Theta\cap\Delta) = \frak{a}(\Theta)\cap\frak{a}(\Delta)$ this follows by
taking perp on both sides and using that $(V+W)^\perp = V^\perp \cap W^\perp$
for $V,W$ linear subspaces.} $\frak{a}_{\Theta \cap \Delta }=\frak{a}%
_{\Theta }+\frak{a}_{\Delta }$. Thus it follows that $Z_{\Theta \cap
\Delta }=Z_{\Theta }\cap Z_{\Delta }$, $K_{\Theta \cap \Delta
}=K_{\Theta }\cap K_{\Delta }$ and $P_{\Theta \cap \Delta
}=P_{\Theta }\cap P_{\Delta }$. For $ H\in \frak{a}$ we denote by
$Z_{H}$, $\mathcal{W}_{H}$ etc.\ the centralizer of $H$,
respectively,  in $G$, $\mathcal{W}$ etc., except when explicitly
noted. When $H\in \mathrm{cl}\frak{a}^{+}$ we put
\[
\Theta (H)=\{\alpha \in \Sigma :\alpha (H)=0\},
\]
and we have $Z_{H}=Z_{\Theta (H)}$, $K_{H}=K_{\Theta (H)}$, $%
N_{H}^{+}=N^{+}(\Theta (H))$ and $\mathcal{W}_{H}=\mathcal{W}_{\Theta (H)}$.

Let $\frak{n}_{\Theta }^{\pm }=\sum_{\alpha \in \Pi ^{\pm }-\langle \Theta
\rangle }\frak{g}_{\alpha }$ and $N_{\Theta }^{\pm }=\exp (\frak{n}_{\Theta
}^{\pm })$. Then $N^{\pm }$ decomposes as $N^{\pm }=N(\Theta )^{\pm
}N_{\Theta }^{\pm }$ where $N(\Theta )^{\pm }$ normalizes $N_{\Theta }^{\pm
} $ and $N(\Theta )^{\pm }\cap N_{\Theta }^{\pm }=1$. We have that $%
\mathfrak{g}=\frak{n}_{\Theta }^{-}\oplus \frak{p}_{\Theta }$, that $%
N_{\Theta }^{-}\cap P_{\Theta }=1$ and also that $P_{\Theta }$ is the
normalizer of $\frak{n}_{\Theta }^{+}$ in $G$. $P_{\Theta }$ decomposes as $%
P_{\Theta }=Z_{\Theta }N_{\Theta }^{+}$, where $Z_{\Theta }$ normalizes $%
N_{\Theta }^{+}$ and $Z_{\Theta }\cap N_{\Theta }^{+}=1$. We write $\frak{p}%
_{\Theta }^{-}=\theta (\frak{p}_{\Theta })$ for the parabolic subalgebra
opposed to $\frak{p}_{\Theta }$. It is conjugate to the parabolic subalgebra
$\frak{p}_{\Theta ^{*}}$ where $\Theta ^{*}=-(w^{-})\Theta $ is the dual to $%
\Theta $ and $w^{-}$ is the principal involution of $\mathcal{W}$. More
precisely, $\frak{p}_{\Theta }^{-}=k\frak{p}_{\Theta ^{*}}$ where $k\in
M^{*} $ is a representative of $w^{-}$. If $P_{\Theta }^{-}$ is the
parabolic subgroup associated to $\frak{p}_{\Theta }^{-}$ then $Z_{\Theta
}=P_{\Theta }\cap P_{\Theta }^{-}$ and $P_{\Theta }^{-}=Z_{\Theta }N_{\Theta
}^{-}$, where $Z_{\Theta }$ normalizes $N_{\Theta }^{-}$ and $Z_{\Theta
}\cap N_{\Theta }^{-}=1$.

The flag manifold of type $\Theta $ is the orbit $\Bbb{F}_{\Theta }=\mathrm{%
Ad}(G)\frak{p}_{\Theta }$ with base point $b_{\Theta }=\frak{p}_{\Theta }$,
which identifies with the homogeneous space $G/P_{\Theta }$. Since the
center of $G$ normalizes $\frak{p}_{\Theta }$, the flag manifold depends
only on the Lie algebra $\frak{g}$ of $G$. The empty set $\Theta =\emptyset $
gives the maximal flag manifold $\Bbb{F}=\Bbb{F}_{\emptyset }$ with
basepoint $b=b_{\emptyset }$. If $\Delta \subset \Theta $ then there is a $G$%
-equivariant projection $\Bbb{F}_{\Delta }\to \Bbb{F}_{\Theta }$ given by $%
gb_{\Delta }\mapsto gb_{\Theta }$, $g\in G$.

The above subalgebras of $\frak{g}$, which are defined by the choice of a
Weyl chamber of $\frak{a}$ and a subset of the associated simple roots, can
be defined alternatively by the choice of an element $H \in \frak{a}$ as
follows. First note that the eigenspaces of $\mathrm{ad}(H)$ in $\frak{g}$
are the weight spaces $\frak{g}_\alpha$, and that the centralizer of $H$ in $%
\frak{g}$ is given by $\frak{z}_H = \sum\{ \frak{g}_\alpha:\, \alpha(H) = 0
\}$, where the sum is taken over $\alpha \in \frak{a}^*$. Now define the
negative and positive nilpotent subalgebras of type $H$ given by
\[
\frak{n}^-_H = \sum\{ \frak{g}_\alpha:\, \alpha(H) < 0 \}, \qquad
\frak{n}^+_H = \sum\{ \frak{g}_\alpha:\, \alpha(H) > 0 \},
\]
and the parabolic subalgebra of type $H$ which is given by
\[
\frak{p}_H = \sum\{ \frak{g}_\alpha:\, \alpha(H) \geq 0 \}.
\]
Denote by $N^\pm_H = \exp(\frak{n}^\pm_H)$ and by $P_H$ the
normalizer in $G$ of $\frak{p}_H$. Note that $\frak{n}^\pm_H$,
$\frak{p}_H$, $N^\pm_H$ and $P_H$ are not centralizers of $H$: these
are the only exceptions for the centralizer notation introduced
above. We have clearly that
\[
\frak{g} = \frak{n}^-_H \oplus \frak{z}_H \oplus \frak{n}^+_H\qquad%
\mbox{and}\qquad \frak{p}_H = \frak{z}_H\oplus \frak{n}^+_H.
\]
Define the flag manifold of type $H$ given by the orbit
\[
{\mathbb F}_H = \mbox{Ad}(G)\frak{p}_H.
\]
Now choose a chamber $\frak{a}^+$ of $\frak{a}$ which contains $H$ in its
closure, consider the simple roots $\Sigma$ associated to $\frak{a}^+$ and
take $\Theta(H) \subset \Sigma$. Since a root $\alpha \in \Theta(H)$ if, and
only if, $\alpha|_{\frak{a}_{\Theta(H)}} = 0$, we have that
\[
\frak{z}_H = \frak{z}_{\Theta(H)},\quad \frak{n}^\pm_H = \frak{n}%
^\pm_{\Theta(H)}, \quad \frak{p}_H = \frak{p}_{\Theta(H)}.
\]
So it follows that
\[
{\mathbb F}_H = {\mathbb F}_{\Theta(H)},
\]
and that the isotropy of $G$ in $\frak{p}_H$ is
\[
P_H = P_{\Theta(H)} = K_{\Theta(H)} A N^+ = K_H A N^+,
\]
since $K_{\Theta(H)} = K_H$. We note that we can proceed reciprocally. That
is, if $\frak{a}^+$ and $\Theta$ are given, we can choose an $H \in \mathrm{%
cl}\frak{a}^+$ such that $\Theta(H) = \Theta$ and describe the objects that
depend on $\frak{a}^+$ and $\Theta$ by $H$ (clearly, such an $H$ is not
unique.) We remark that the map
\begin{equation}  \label{eq:mergulho-em-s}
{\mathbb F}_H \to \mathfrak{s},\qquad k \frak{p}_H \mapsto \mbox{Ad}(k)
H,\quad\mbox{where } k \in K,
\end{equation}
gives an embeeding of ${\mathbb F}_H$ in $\mathfrak{s}$ (see Proposition 2.1
of \cite{dkv}). In fact, the isotropy of $K$ at $H$ is $K_H = K_{\Theta(H)}$
which is, by the above comments, the isotropy of $K$ at $\frak{p}_H$. Define
the negative parabolic subalgebra of type $H$ by
\[
\frak{p}^-_H = \sum\{ \frak{g}_\alpha:\, \alpha(H) \leq 0 \}
\]
and denote by $P^-_H$ its normalizer in $G$. Then we have that $P^-_H =
P^-_{\Theta(H)}$.

\section{Gradient flows on flag manifolds\label{seclineariz}}

The main result of this section is the construction in Theorem \ref
{teolineariz} of equivariant linearizations for the gradient flows on the
flag manifolds induced by the split elements in $\frak{g}$. In the next
section we plug this linearization on the fibers of the flag bundles to get
our main tool in the determination of the Conley indices.

\subsection{Fixed points, stable and unstable sets}

An split element $H\in \mathrm{cl}\frak{a}^{+}$ induces a vector field $%
\widetilde{H}$ on a flag manifold ${\mathbb F}_{\Theta }$ with flow $\exp
(tH)$. This is a gradient vector field with respect to a given Riemannian
metric on ${\mathbb F}_{\Theta }$ (see \cite{dkv}, Section 3). The connected
sets of fixed point of this flow are given by
\[
\mathrm{fix}_{\Theta }(H,w)=Z_{H} wb_{\Theta }=K_{H} wb_{\Theta },
\]
so that they are in bijection with the cosets in $\mathcal{W}_{H}\backslash
\mathcal{W}/\mathcal{W}_{\Theta }$. Each $w$-fixed point connected set has
stable manifold given by
\[
\mathrm{st}_{\Theta }(H,w)=N_{H}^{-} \mathrm{fix}_{\Theta }(H,w)=P_{H}^{-}
wb_{\Theta },
\]
whose union gives the Bruhat decomposition of ${\mathbb F}_{\Theta }$:
\[
{\mathbb F}_{\Theta }=\coprod_{\mathcal{W}_{H}\backslash \mathcal{W}/%
\mathcal{W}_{\Theta }}\mathrm{st}_{\Theta }(H,w)=\coprod_{\mathcal{W}%
_{H}\backslash \mathcal{W}/\mathcal{W}_{\Theta }}P_{H}^{-} wb_{\Theta }.
\]
The unstable manifold is
\[
\mathrm{un}_{\Theta }(H,w)=N_{H}^{+} \mathrm{fix}_{\Theta }(H,w)=P_{H}
wb_{\Theta }.
\]
Since the centralizer $Z_{H}$ of $H$ leaves $\mathrm{fix}_{\Theta }(H,w)$
invariant and normalizes both $N_{H}^{-}$ and $N_{H}^{+}$, it follows $%
\mathrm{st}_{\Theta }(H,w)$ and $\mathrm{un}_{\Theta }(H,w)$ are $Z_{H}$%
-invariant. We note that these fixed points and (un)stable sets remain the
same if $H$ is replaced by some $H^{\prime }\in \mathrm{cl}\frak{a}^{+}$
such that $\Theta (H^{\prime}) = \Theta(H)$.

For $X\in \frak{g}$ we denote with $\widetilde{X}$ the vector field induced
by $X$ on ${\mathbb F}_{\Theta }$, and for a subset $\frak{l}\subset \frak{g}
$ and $x\in {\mathbb F}_{\Theta }$ we put
\[
\frak{l}\cdot x=\{\widetilde{X}\left( x\right) \in T_{x}{\mathbb F}_{\Theta
}:X\in \frak{l}\}.
\]
With this notation, for some $x\in \mathrm{fix}_{\Theta }(H,w)$, the
tangent space to $\mathrm{fix}_{\Theta }(H,w)$ at $x$ is given by
$\frak{z}_{H}\cdot x$, while the tangent spaces at $x$ of the
unstable and stable manifolds is given, respectively, by
$\frak{n}^\pm_H\cdot x$. These tangent spaces form the stable and
unstable vector bundles
\[
V_{\Theta }^{+}(H,w)\rightarrow \mathrm{fix}_{\Theta }(H,w)\qquad \mathrm{and%
}\qquad V_{\Theta }^{-}(H,w)\rightarrow \mathrm{fix}_{\Theta }(H,w)
\]
with fibers $\frak{n}^\pm_H\cdot x$, $x\in \mathrm{fix}_{\Theta
}(H,w)$, respectively. Their Whitney sum
\begin{equation}  \label{eqwhitneysum}
V_{\Theta }\left( H,w\right) =V_{\Theta }^{+}(H,w)\oplus V_{\Theta }^{-}(H,w)
\end{equation}
is a normal bundle of $\mathrm{fix}_{\Theta }(H,w)$ in ${\mathbb
F}_{\Theta }$. Note that thanks to the well known translation formula for
vector fields
\[
dg_{x}(\widetilde{X}(x))=( \mathrm{Ad}(g)X)^{\sim }(gx),
\]
for $g\in G$, it follows that the bundles $V_{\Theta }^{\pm }(H,w)$ are $%
Z_{H}$-invariant.

\subsection{Linearization\label{subseclinflag}}

In the sequel a special role is played by the subgroup $L_{H}$ of $Z_{H}$, $%
H\in \mathrm{cl}\frak{a}^{+}$, defined by
\begin{equation}
L_{H}=K_{H}A_{\Theta (H)}\subset Z_{H}.  \label{forconforme}
\end{equation}
Note that $L_{H}$ is the direct product of a compact group $K_{H}$
and a vector group $A_{\Theta (H)}$, hence we call it the conformal
part of $Z_{H}$. In case $H$ is regular $L_{H}=Z_{H}=MA$. Also, by
the structure of the parabolic subgroups we have
\[
{\mathrm{fix}_{\Theta }}(H,w)=L_{H}wb_{\Theta }\qquad
\mathrm{and}\qquad \mathrm{st}_{\Theta }(H,w)=L_{H}N^{-}_Hwb_{\Theta
}.
\]

Now we start the main construction of this section, namely of a
mapping
conjugating the flow $\exp (tH)$ around a fixed point set $\mathrm{fix}%
_{\Theta }(H,w)$ with a linear flow on the above normal bundle of $\mathrm{fix}%
_{\Theta }(H,w)$. The existence of such a conjugation is a
consequence of general theory on Morse-Bott theory (see \cite{dkv},
Corollary 1.5). However we itend to carry this construction to flag
bundles, so we require the conjugation map to be equivariant with
respect to $L_{H}$. Unfortunately our
method does not yield, in general, equivariance with respect to the whole $%
Z_{H}$, as desirable (see the next remark). Below we discuss some cases where $%
Z_{H}$-equivariance holds.

The strategy to get a conjugation consists in defining for each
$x\in \mathrm{fix}_{\Theta }(H,w)$ a subspace
\[
\frak{l}_{x}\subset \frak{n}^-_H\oplus \frak{n}^+_H
\]
satisfying the following requirements:

\begin{enumerate}
\item  $\frak{l}_{x}\cdot x=V_{\Theta }\left( H,w\right)
_{x}=\left( \frak{n}^-_H\oplus \frak{n}^+_H\right) \cdot x$.

\item  The map $X\in \frak{l}_{x}\mapsto \widetilde{X}\left( x\right) \in
V_{\Theta }\left( H,w\right) $ is injective, and hence a bijection.

\item  For every $g\in L_{H}$, $\mathrm{Ad}\left( g\right) \frak{l}_{x}=%
\frak{l}_{gx}$. (This is to ensure equivariance.)
\end{enumerate}
Once we have the subspaces $\frak{l}_{x}$ it is easy to get a
conjugation by putting
\[
\psi :V_{\Theta }(H,w)\to {\mathbb F}_{\Theta }, \qquad \psi \left(
v\right) =\exp (X)x
\]
where $X\in \frak{l}_{x}$ is such that $v=\widetilde{X}\left(
x\right) $. By the third condition and the translation formula for
vector fields, it follows that $\psi $ is $L_{H}$-equivariant, that
is
\[
\psi \circ dg=g\circ \psi, \qquad g\in L_{H}.
\]

We define the subspaces ${\frak l}_x$, $x \in {{\rm
fix}_\Theta(H,w)}$, by
$$
{\frak l}_x = \mbox{Ad}(g){\frak l}_{w b_\Theta}, \qquad {\frak
l}_{w b_\Theta} = w{\frak n}^-_\Theta \cap \left( \frak{n}^-_H\oplus
\frak{n}^+_H\right),
$$
where $x = g w b_\Theta$, with $g \in L_H$.  In order to show that
${\frak l}_x$ is well defined, we only have to check that the
isotropy of $w b_\Theta$ in $L_H$ normalizes $w{\frak n}^-_\Theta$,
since $L_H$ already normalizes $ \frak{n}^-_H\oplus \frak{n}^+_H$.
\begin{lema}
The isotropy of $w b_\Theta$ in $L_H$ normalizes $w{\frak
n}^-_\Theta$.
\end{lema}
\begin{profe}
Clearly, this is isotropy is given by $P^w_\Theta \cap L_H$, where
the superscript $w$ denotes conjugation by some representative of
$w$ in $M^*$. Using the Iwasawa decompositions
\[
P^w_\Theta = K_\Theta^w A N^w, \qquad Z^w_\Theta = K_\Theta^w A
N(\Theta)^w
\]
and the uniqueness of the Iwasawa decomposition $G = KAN^w$, we have
that
$$
P^w_\Theta \cap L_H = (K_\Theta^w \cap K_H) A_{\Theta(H)} \subset
Z_\Theta^w,
$$
which normalizes $w{\frak n}^-_\Theta$.
\end{profe}

In order to verify conditions (1) and (2) stated above, we only need
to check that the isotropy subalgebra of ${\frak g}$ at $x$ is
complemented by ${\frak l}_x$. In fact, the isotropy subalgebra at
$x = g wb_\Theta$, $g \in L_H$, is given by $\mbox{Ad}(g)w{\frak
p}_\Theta$. It is clearly complemented in ${\frak g}$ by ${\frak
l}_x$, since ${\frak g} = {\frak p}_\Theta \oplus {\frak
n}^-_\Theta$. Condition (3) is immediate from the definition of
${\frak l}_x$.

Now we can state the conjugation result.

\begin{teorema}
\label{teolineariz}The map $\psi :V_{\Theta }(H,w)\to {\mathbb
F}_{\Theta }$ is $L_{H}$-equivariant local diffeomorphism and
satisfies:

\begin{itemize}
\item[i)]  Its restriction to a neighborhood of the zero section is a
diffeomorphism onto a neighborhood of ${\mathrm{fix}_{\Theta
}}(H,w)$.

\item[ii)]  Its restriction to $V_{\Theta }^{-}(H,w)$ and $V_{\Theta
}^{+}(H,w)$ are diffeomorphisms onto $\mathrm{st}_{\Theta }(H,w)$ and $%
\mathrm{un}_{\Theta }(H,w)$, respectively.
\end{itemize}
\end{teorema}

\begin{profe}
The $L_{H}$-equivariance follows by construction and it is also
immediate that $\psi$ is the identity on the null section. Next we
show that $\psi$ is a local diffeomorphism. By the
$L_H$-equivariance and the inverse function theorem, it is enough to
show that the differential ${\rm d}_v\psi$ is an isomorphism for $v
= {\widetilde{X}(wb_\Theta)}$, where $X \in {\frak l}_{wb_\Theta}$.
We first note that the dimensions of the manifolds $V_\Theta(H,w)$
and ${\mathbb F}_\Theta$ are the same, since $V_\Theta(H,w)$ is a
normal bundle. Hence we only need to show that the differential is
surjective. Taking $Z \in {\frak l}^0_{wb_\Theta} = w{\frak
n}^-_\Theta \cap {\frak z}_H$, we have that
$$u(t) = \widetilde{X}(\exp(tZ)wb_\Theta)$$
is a curve that is transversal to the fibers of $V_\Theta(H,w)$.
Thus we can compute
$$
\left. \frac{\rm{d}}{\rm{d}t} \right|_{t=0}  \psi(u(t)) = \left.
\frac{\rm{d}}{\rm{d}t} \right|_{t=0} \exp(X)\exp(tZ)wb_\Theta = {\rm
d}_{wb_\Theta}\exp(X)\, \widetilde{Z}(wb_\Theta).
$$
On the other hand, taking $Y \in {\frak l}_{wb_\Theta}$, we have
that
$$v(t) = (\log(\exp(X)\exp(tY))^\sim(wb_\Theta)$$
is a curve of $V_\Theta(H,w)$ in the fiber over $wb_\Theta$.  Thus
we can compute
$$
\left. \frac{\rm{d}}{\rm{d}t} \right|_{t=0}  \psi(v(t)) =  \left.
\frac{\rm{d}}{\rm{d}t} \right|_{t=0} \exp(X)\exp(tY) wb_\Theta =
{\rm d}_{wb_\Theta}\exp(X)\, \widetilde{Y}(wb_\Theta).
$$
Since ${\rm d}_{wb_\Theta}\exp(X)$ is an isomorphism, the dimension
of the image of ${\rm d}_v\psi$, where $v =
{\widetilde{X}(wb_\Theta)}$, is greater than the dimension of the
space spanned at $wb_\Theta$ by the induced vectors fields of
${\frak l}_{wb_\Theta} \oplus {\frak l}^0_{wb_\Theta}$. Since this
is precisely the dimension of ${\mathbb F}_\Theta$ it follows that
${\rm d}_v\psi$ is surjective.

For item (i), suppose by contradiction that there is no neighborhood
of ${{\rm fix}_\Theta(H,w)}$ in $V_\Theta(H,w)$ in which $\psi$ is
injective. Then we have sequences $u_k, v_k \in V_\Theta(H,w)$ such
that $u_k \neq v_k$, $\psi(u_k) = \psi(v_k)$ and $u_k,v_k \to {{\rm
fix}_\Theta(H,w)}$ when $k\to \infty$. Since ${{\rm
fix}_\Theta(H,w)}$ is compact, we can take a compact neighborhood of
the null section in $V_\Theta(H,w)$. Thus we have, taking
subsequences, that $u_k \to x$, $v_k \to y$, where $x,y \in {{\rm
fix}_\Theta(H,w)}$. But then
$$x = \psi(x) = \lim_k \psi(u_k) = \lim_k \psi(v_k) = \psi(y) = y.$$
Since $u_k \neq v_k$ the map $\psi$ fails to be locally injective
around $x=y$, contradicting the local diffeomorphism property
already proved.

For item (ii) we first compute the image of $\psi$ restricted to
$V_\Theta^-(H,w)$. Note that
\[
V_\Theta^-(H,w)_{wb_\Theta} = \frak{n}^-_H \cdot wb_\Theta
\]
and thus
\[
V_\Theta^-(H,w) = L_H \cdot V_\Theta^-(H,w)_{wb_\Theta} = L_H \cdot
\left( \frak{n}^-_H \cdot wb_\Theta \right).
\]
Therefore, by the $L_H$-equivariance, we have that
\[
\psi (V_{\Theta }^{-}(H,w)) = L_{H}N_{H}^{-}\cdot wb_{\Theta } =
N_{H}^{-}\cdot \mathrm{fix}_{\Theta }(H,w) = \mathrm{st}_{\Theta
}(H,w),
\]
so that the map is onto. For the injectivity suppose that
\[
g\exp (X)\cdot wb_{\Theta }=k\exp (Y)\cdot wb_{\Theta }
\]
where $g,k\in L_{H}$ and $X,Y\in w\frak{n}_{\Theta}^{-}\cap
\frak{n}^-_H$. Applying $\exp tH$ to both sides and making
$t\rightarrow +\infty $, the left hand side converges to $g\cdot
wb_{\Theta }$ and the right hand side to $k\cdot wb_{\Theta }$.
Hence, $g\cdot wb_{\Theta }=k\cdot wb_{\Theta }$. Hence $g^{-1}k\in
wP_{\Theta }w^{-1}$. It follows that
\[
\exp (-X)\exp (\mathrm{Ad}(g^{-1}k)Y)\in w(P_{\Theta }\cap N_{\Theta
}^{-})w^{-1}=1.
\]
Since $\exp$ restricted to ${\frak n}^-$ is a diffeomorphism, it
follows that $\mathrm{Ad}(g)X=\mathrm{Ad}(h)Y$ and
$(\mathrm{Ad}(g)X)^{\sim }(g\cdot wb_{\Theta
})=(\mathrm{Ad}(h)Y)^{\sim }(k\cdot wb_{\Theta })$. That is, the map
is injective. The proof for $V_{\Theta }^{+}(H,w)$ is analogous.
\end{profe}

\begin{remarque}
The bundles $V_{\Theta }^{\pm }(H,w)$ as well as the stable and unstable
sets $\mathrm{st}_{\Theta }(H,w)$ and $\mathrm{un}_{\Theta }(H,w)$ are $%
Z_{H} $-invariant. Despite this, the conjugation $\psi $ is not in general $%
Z_{H}$-equivariant because it may not be true that the subspace
\[
\frak{l}_{wb_{\Theta }}=w\frak{n}_{\Theta }^{-}\cap
(\frak{n}^-_H\oplus \frak{n}^+_H)
\]
is invariant under the isotropy at $wb_{\Theta }$ of the
$Z_{H}$-action. An example of the non invariance is given next in
the maximal flag manifold ($\Theta = \emptyset$).
\end{remarque}

\begin{exemplo}
The isotropy at $wb$ of the $Z_{H}$-action is $Z_{H}\cap P^{w}$ whose Lie
algebra is
\[
\frak{z}_{H}\cap w\frak{p}=
w\left( {\textstyle \sum_{\beta }} \frak{g}_{\beta }\right)
\]
with the sum extended to the roots $\beta \geq 0$ such that $\beta
(w^{-1}H)=0$. On the other hand
\[
\frak{l}_{wb}=w\frak{n}^{-}\cap (\frak{n}^-_H\oplus \frak{n}^+_H)=
w\left( {\textstyle \sum_{\alpha }} \frak{g}_{\alpha}\right)
\]
with the sum over $\alpha <0$ with $\alpha (w^{-1}H)\neq 0$. Clearly $\frak{l%
}_{wb}$ is normalized by $\frak{z}_{H}\cap w\frak{p}$ if it is invariant by $%
Z_{H}\cap P^{w}$.

Now let $X_{\alpha }\in \frak{g}_{\alpha }$, where $\alpha <0,\,\alpha
(w^{-1}H)\neq 0$, and $X_{\beta }\in \frak{g}_{\beta }$, where $\beta \geq
0,\,\beta (w^{-1}H)=0$, be such that $[X_{\alpha },X_{\beta }]\neq 0$. If $%
\frak{z}_{H}\cap w\frak{p}$ normalizes $\frak{l}_{wb}$ then
\[
w[X_{\alpha },X_{\beta }]=[wX_{\alpha },wX_{\beta }]\in \frak{l}_{wb}
\]
which happens if, and only if,
\[
(\alpha +\beta )(w^{-1}H)\neq 0\qquad \mbox{and}\qquad \alpha +\beta <0,
\]
since $[X_{\alpha },X_{\beta }]\in \frak{g}_{\alpha +\beta }$. Keeping this
in mind, our example is given for $\frak{g}=\frak{sl}(3,\mathbb{R})$ with
the usual choices of chambers and roots. Let $w\in \mathcal{W}$ be the
permutation $w=\left( 23\right) =w^{-1}$, and
\[
H=\left(
\begin{array}{ccc}
-1 &  &  \\
& -1 &  \\
&  & 2
\end{array}
\right) ,\quad X_{\alpha }=\left(
\begin{array}{ccc}
0 & 0 & 0 \\
1 & 0 & 0 \\
0 & 0 & 0
\end{array}
\right) ,\quad X_{\beta }=\left(
\begin{array}{ccc}
0 & 0 & 1 \\
0 & 0 & 0 \\
0 & 0 & 0
\end{array}
\right) ,
\]
where $\alpha =\alpha _{21}=-\alpha _{12}$, $\beta =\alpha _{13}=\alpha
_{12}+\alpha _{23}$. Thus $\alpha +\beta =\alpha _{23}>0$, $\alpha <0$ with $%
\alpha (w^{-1}H)=3\neq 0$, $\beta >0$ with $\beta (w^{-1}H)=0$ and $%
[X_{\alpha },X_{\beta }]\neq 0$. By the general remarks above it follows
that $\frak{z}_{H}\cap w\frak{p}$ does not normalize $\frak{l}_{wb}$ for
this choice of $H$ and $w$.
\end{exemplo}

In conclusion we mention some cases where the conjugation is $Z_{H}$
invariant.

\begin{corolario}
\label{prop2lineariz}If $H$ is regular then $Z_{H}=L_{H}$ so that $\psi $ is
$Z_{H}$-equivariant.
\end{corolario}

Another case where $Z_{H}$-invariance holds is for the attractor fixed point
component $\mathrm{fix}_{\Theta }\left( H,1\right) $ when the flag manifold $%
\Bbb{F}_{\Theta }$ either projects onto $\Bbb{F}_{\Theta \left( H\right) }$
or $\Bbb{F}_{\Theta \left( H\right) }$ projects onto $\Bbb{F}_{\Theta }$.

\begin{corolario}
\label{prop1lineariz}Let $w=1$ and suppose that
\[
\Theta \subset \Theta (H)\qquad \mathrm{or}\qquad \Theta (H)\subset \Theta .
\]
Then $\psi :V_{\Theta }(H,1)=V_{\Theta }(H,1)^{-}\rightarrow \Bbb{F}_{\Theta
}$ is $Z_{H}$-equivariant and is a global diffeomorphism onto the open and
dense subset $\mathrm{st}_{\Theta }(H,1)$.
\end{corolario}

\begin{profe}
To prove $Z_{H}$-invariance it is enough to check that the subalgebra
\[
\frak{n}_{\Theta }^{-}\cap \frak{n}^-_H
\]
is invariant under $\mathrm{Ad}\left( Z_{H}\cap P_{\Theta }\right)
$. In general $\frak{n}^-_H$ is $\mathrm{Ad}\left( Z_{H}\right)
$-invariant.
Now, if $\Theta (H)\subset \Theta $ then $Z_{H}\subset Z_{\Theta }$, hence $%
Z_{H}\cap P_{\Theta }\subset Z_{\Theta }$ which normalizes $\frak{n}_{\Theta
}^{-}$, so that $\frak{n}_{\Theta }^{-}\cap \frak{n}^-_H$ is $\mathrm{Ad}%
\left( Z_{H}\cap P_{\Theta }\right) $-invariant. On the other hand when $%
\Theta \subset \Theta (H)$ then $\frak{n}^-_H\subset
\frak{n}_{\Theta }^{-}$ so that $\frak{n}_{\Theta }^{-}\cap
\frak{n}^-_H=\frak{n}^-_H$, which is normalized by $Z_{H}$ showing
the result in this case as well. The
proof then follows closely the proof of Theorem \ref{teolineariz} taking $%
w=1 $ and replacing $L_{H}$ by $Z_{H}$.
\end{profe}

\begin{remarque}
A linearization similar to the above one was already given in Proposition
3.6 of \cite{dkv}. However that linearization is not suitable for our
purposes for two reasons. First it does not cover (in the nonregular case)
the whole fixed point component, but only a dense subset of it. Also, it is $%
wZ_{H}w^{-1}$-equivariant instead of $Z_{H}$-equivariant as we will need
latter.
\end{remarque}

\subsection{Fixed point components as flag manifolds}

It is an interesting fact that each fixed point component ${\mathrm{fix}%
_{\Theta }}(H,w)$ in $\Bbb{F}_{\Theta }$ is actually diffeomorphic to a flag
manifold of a certain subgroup of $G$. In fact, for $\Delta \subset \Sigma$
let $\frak{g}(\Delta)$ be the semi-simple subalgebra of $\frak{g}$ of type $%
\Delta$. For $H_0 \in \frak{a}(\Delta)$ let $\mathbb{F}(\Delta)_{H_0}$ be
the flag manifold of $\frak{g}(\Delta)$ of type ${H_0}$. Note that we do not
need to worry in which chamber of $\frak{a}(\Delta)$ the element $H_0$ lies.
Denote by $\pi_\Delta: \frak{a} \to \frak{a}(\Delta)$ the orthogonal
projection.

\begin{proposicao}
\label{prop:pto-fixo-flag} Let $H\in \mathrm{cl}\frak{a}^{+}$, $w\in W$. Put
$\Delta =\Theta (\phi )$ and take $H_{\Theta }\in \mathrm{cl}\frak{a}^{+}$
such that $\Theta (H_{\Theta })=\Theta $. Then we have a diffeomorphism
\[
\mathrm{fix}(H,w)_{\Theta }\simeq \mathbb{F}(\Delta)_{H_{0}},
\]
where $H_{0}$ is the orthogonal projection of $wH_{\Theta }$ in $\frak{a}%
\left( \Delta \right) $. Furthermore, in the maximal flag manifold $%
\mathbb{F}$, there is a $K_{H}$-equivariant diffeomorphism between any two $%
\mathrm{fix}(H,w)$, $w\in W$, which are thus diffeomorphic to the maximal
flag manifold of $\frak{g}(\Delta )$.
\end{proposicao}

\begin{profe}
Decompose $wH_{\Theta }=H_{0}+H_{1}$, where $H_{0}\in \frak{a}(\Delta )$, $%
H_{1}\in \frak{a}_{\Delta }$. Using that $K_{H}=K(\Delta )M$ we have that
\[
{\mathrm{fix}_{\Theta }(H,w)}=K_{H}wb_{\Theta }=K(\Delta )wb_{\Theta }.
\]
Using that $K_{\Theta }=K_{H_{\Theta }}$ it follows that the isotropy of $K$
at $wb_{\Theta }$ is $K\cap wP_{\Theta }{w^{-1}}=wK_{\Theta }{w^{-1}}%
=K_{wH_{\Theta }}$. So the isotropy of $K(\Delta )$ at $wb_{\Theta }$ is $%
K(\Delta )_{wH_{\Theta }}=K(\Delta )_{H_{0}}$, where the last equality
follows from $K(\Delta )\subset K_{\Delta }$ so it centralizes $H_{1}\in
\frak{a}_{\Delta }$. It follows that
\[
{\mathrm{fix}_{\Theta }(H,w)}\to \mathfrak{s}(\Delta),\qquad kwb_{\Theta
}\mapsto \mbox{Ad}(k)H_{0},\quad \mbox{where }k\in K(\Delta ),
\]
is an embeeding. Equation (\ref{eq:mergulho-em-s}) applied to $\frak{g}%
(\Delta )$ gives that the image of this embeeding is diffeomorphic to the
flag manifold ${\mathbb F}(\Delta )_{H_{0}}$. This proves the first
assertion. For the second assertion, we note that
\[
{\mathbb F}\to {\mathbb F},\quad kb\mapsto kwb,\quad \mbox{where }k\in K,
\]
is a well defined $K$-invariant diffeomorphism of ${\mathbb F}$. In fact, if
$kb=k^{\prime }b$ then $k^{\prime }=km$, $m\in M$. Since $w$ normalizes $M$,
it follows that $k^{\prime }wb=kwb$. The above map clearly restricts to a $%
K_{H}$-invariant diffeomorphism between $\mathrm{fix}(H,1)=K_{H}b$ and $%
\mathrm{fix}(H,w)=K_{H}wb$. By the first part of the proof, we have that $%
\mathrm{fix}(H,1)\simeq \mathbb{F}(\Delta)_{\pi _{\Delta }(H_{\Theta })}$,
where $H_{\Theta }$ is regular since we are in the maximal flag manifold. We
claim that $\pi _{\Delta }(H_{\Theta })=H_{0}$ is regular in $\frak{a}%
(\Delta )$, which implies that $\mathbb{F}(\Delta)_{H_{0}}$ is the maximal
flag manifold of $\frak{g}(\Delta )$. In fact, for all $\alpha \in \langle
\Delta \rangle $ we have $\alpha |_{\frak{a}_{\Delta }}=0$ so that $\alpha
(H_{0})=\alpha (H_{\Theta })\neq 0$.
\end{profe}

\subsection{Schubert cells\label{subsecschub}}

We recall some facts about the closure of the stable sets $\mathrm{st}%
_{\Theta }(H,w)$. If $H$ is regular then the closure $\mathrm{cl}(\mathrm{st}%
_{\Theta }(H,w))$ is known as a Schubert cell in the flag manifold $\Bbb{F}%
_{\Theta }$. In this case a result that goes back to Borel-Tits \cite{botit}
ensures that a Schubert cell is a union of Bruhat cells. Precisely,
\begin{equation}
\mathrm{cl}(\mathrm{st}_{\Theta }(H,w))=\bigcup_{s\geq w}\mathrm{st}_{\Theta
}(H,s)=\bigcup_{s\geq w}N^{-}\left( s\cdot b_{\Theta }\right)
\label{forschubert}
\end{equation}
where the order in $\mathcal{W}$ is Bruhat-Chevalley. For nonregular $H$ the
same equality is true because a fixed point component $\mathrm{fix}_{\Theta
}\left( H,w\right) $ is a $K_{H}-$orbit and
\[
\mathrm{st}_{\Theta }\left( H,w\right) =N_{H}^{-}K_{H}wb_{\Theta
}=K_{H}N^{-}wb_{\Theta },
\]
that is, $\mathrm{st}_{\Theta }\left( H,w\right) $ is the union of $K_{H}$%
-translations of Bruhat cells. Hence (\ref{forschubert}) holds as well after
taking closures, since $K_H$ is compact.

These descriptions of Schubert cells are group theoretic. On the other hand
on the maximal flag manifold there is the following alternative description
of \cite{sm:order}: Fix a simple system of roots $\Sigma $, and for a finite
sequence $\alpha _{1},\ldots ,\alpha _{n}$ in $\Sigma $ we let $s_{1},\ldots
,s_{n}$ be the reflections with respect to these roots. Then we write $\Bbb{F%
}_{i}=\Bbb{F}_{\{\alpha _{i}\}}$ and put $\pi _{i}:\Bbb{F}\rightarrow \Bbb{F}%
_{i}$ for the canonical projection. Accordingly, we write $\gamma _{i}=\pi
_{i}^{-1}\pi _{i}$ for the map that to a subset $A\subset \Bbb{F}$
associates the union of the fibers crossing $A$. Then it follows by the
results of \cite{sm:order} that
\[
\mathrm{cl}(\mathrm{st}(H,w))=\gamma _{1}\cdots \gamma _{n}\left( \mathrm{fix%
}\left( H,w^{-}\right) \right) ,
\]
where $\gamma _{1},\ldots ,\gamma _{n}$ is taken from a reduced expression $%
w^{-}w=s_{n}\cdots s_{1}$.

\section{Morse decomposition on flag bundles}

From now on we fix the setup stated in the introduction: A principal bundle $%
\pi :Q\to X$ with semi-simple structural group $G$ and compact Hausdorff
base space $X$. A right invariant flow $\phi _{t}$ on $Q$ which is chain
transitive on $X$. The flow $\phi _{t}$ induces a flow in the flag bundle $%
\Bbb{E}_{\Theta }$. The following result, proved in \cite{smbflow} and \cite
{msm}, gives the chain components and Morse decomposition of this induced
flow.

\begin{teorema}
\label{teo:tipo-parabol-flow}The flow $\phi _{t}$ on the flag bundle ${\Bbb{E%
}_{\Theta }}\to X$ admits a finest Morse decomposition whose Morse
components are determined as follows.

\begin{itemize}
\item[i)]  There exists $H_{\phi }\in \mathrm{cl}\frak{a}^{+}$ and a $\phi
_{t}$-invariant $Z_{H_{\phi }}$-reduction $Q_{\phi }$ of $Q$ such that the
Morse components are parameterized by the Weyl group $\mathcal{W}$ and each
component $\mathcal{M}_{\Theta }^{w}$ is given by
\[
{\mathcal{M}}_{\Theta }^{w}=Q_{\phi }\cdot {\mathrm{fix}}(H_{\phi
},w)_{\Theta }=Q_{\phi }\cdot wb_{\Theta }.
\]

\item[ii)]  There is just one attractor component $\mathcal{M}_{\Theta }^{+}=%
\mathcal{M}_{\Theta }^{1}$ and only one repeller $\mathcal{M}_{\Theta }^{-}=%
\mathcal{M}_{\Theta }^{w^{-}}$.

\item[iii)]  The dynamical ordering of the Morse components is given by the
algebraic Bruhat-Chevalley order of $\mathcal{W}$.

\item[iv)]  There exists a $K_{H_{\phi }}$-reduction $R_{\phi }$ of $Q_{\phi
}$ (which is not necessarily $\phi _{t}$-invariant) such that each Morse
component is given by the orbit
\[
{\mathcal{M}}_{\Theta }^{w}=R_{\phi }\cdot {\mathrm{fix}}(H_{\phi
},w)_{\Theta }=R_{\phi }\cdot wb_{\Theta }.
\]
\end{itemize}
\end{teorema}

We say that $H_{\phi }$ a characteristic element for the flow. The set of
simple roots $\Theta (\phi )=\Theta (H)$ is called the parabolic type of the
flow and the $\phi _{t}$-invariant $Z_{H}$-reduction $Q_{\phi }$ is called a
block reduction of the flow (cf. \cite{ll1}). When the characteristic
element $H_{\phi }$ is regular we say that the flow $\phi _{t}$ is regular
as well.

In the sequel we say that a subset $A\subset {\mathbb E}_{\Theta }$ is a
section over $X$ if it is the image of a continuous section of the bundle ${%
\mathbb
E}_{\Theta }\to X$. If a subset is a section then it meets each fiber in a
singleton. For example the $\mathcal{M}_{\Theta \left( \phi \right) }^{+}$
is a section of the flag bundle ${\mathbb
E}_{\Theta \left( \phi \right) }\to X$ associated to the parabolic type of $%
\phi $.

Next we show the interesting fact that the Morse components $\mathcal{M}%
_{\Theta }^{w}$ can be viewed as a flag subbundles of ${\mathbb E}_{\Theta }$
(cf. Proposition \ref{prop:pto-fixo-flag} above).

\begin{proposicao}
\label{prop:compon-subflags} Put $\Delta =\Theta (\phi )$ and take $%
H_{\Theta }\in \mathrm{cl}\frak{a}^{+}$ such that $\Theta (H_{\Theta
})=\Theta $. Then the following statements hold.

\begin{description}
\item[i)]  A Morse component $\mathcal{M}_{\Theta }^{w}$ is an associated
bundle of $Q_{\phi }\to X$ with typical fiber $\mathbb{F}(\Delta)_{H_{0}}$
where $H_{0}$ is the orthogonal projection of $wH_{\Theta }$ in $\frak{a}%
\left( \Delta \right) $.

\item[ii)]  In the maximal flag bundle ${\mathbb E}$, there is a bundle
homeomorphism between any two Morse components $\mathcal{M}^{w}$, $w\in
\mathcal{W}$. The typical fiber of any one of these bundles is the maximal
flag manifold ${\mathbb F}(\Delta )$ of $G\left( \Delta \right) $.

\item[iii)]  If the flow is regular, then each Morse component $\mathcal{M}%
_{\Theta }^{w}$ is a section over $X$.
\end{description}
\end{proposicao}

\begin{profe}
By Theorem \ref{teo:tipo-parabol-flow}, we have $\mathcal{M}_{\Theta
}^{w}=Q_{\phi }\cdot wb_{\Theta }$ hence $\mathcal{M}^{w}\to X$ is a fiber
bundle associated to $Q_{\phi }\to X$ with typical fiber $Z_{H}wb_{\Theta }=%
\mathrm{fix}_{\Theta }(H,w)$. Using that $\Delta =\Theta (\phi )=\Theta
(H_{\phi })$ the result follows from the first part of Proposition \ref
{prop:pto-fixo-flag}. In order to construct a fiber-bundle homeomorphism
between the components $\mathcal{M}^{w}\subset {\mathbb E}$ we consider the $%
K_{H}$ reduction $R_{\phi }$ of $Q_{\phi }$ given by item (iv) of Theorem
\ref{teo:tipo-parabol-flow} and also the $K_{H}$-equivariant diffeomorphism
between any two $\mathrm{fix}(H,w)$ given by Proposition \ref
{prop:pto-fixo-flag}. Plugging this diffeomorphism on the fibers of the flag
bundles through the $R_{\phi }$ reduction yields the desired bundle
homeomorphism. Finally in the regular case the typical fibers are just a
points.
\end{profe}

\section{Stable and unstable sets\label{secstable}}

The stable and unstable sets
\[
\mathrm{st}(\mathcal{M})=\{x\in E:\omega (x)\subset \mathcal{M}\}\qquad
\mathrm{un}(\mathcal{M})=\{x\in E:\omega ^{*}(x)\subset \mathcal{M}\}
\]
of a Morse component $\mathcal{M}$ in a flag bundle are given in a similar
fashion as the components themselves. Namely, by plugging fiberwise stable
and unstable sets on the flag manifolds. We work out the stable sets. The
unstable ones are obtained by symmetry.

We start by recalling a result of \cite{smbflow}. Following the notation of
Subsection \ref{subsecschub} let $s_{1},\ldots ,s_{n}$ be reflections with
respect to the simple roots $\alpha _{1},\ldots ,\alpha _{n}$ and write ${%
\mathbb E}_{i}={\mathbb E}_{\{\alpha _{i}\}}$ for the corresponding flag
budles. Put put $\pi _{i}:{\mathbb E}\rightarrow {\mathbb E}_{i}$ for the
canonical projection and $\gamma _{i}=\pi _{i}^{-1}\pi _{i}$ for the map
that exhausts fibers. Note that each map $\gamma _{i}$ preserves the fibers
of ${\mathbb E}\rightarrow X$.

The following statement is Proposition 9.9 of \cite{smbflow}. It
characterizes the domains of attraction in the maximal flag bundle.

\begin{proposicao}
\label{propdmatalg}The domain of attraction of $\mathcal{M}^{w}$ is given by
\begin{equation}
\mathcal{A}\left( \mathcal{M}^{w}\right) =\gamma _{1}\cdots \gamma
_{n}\left( \mathcal{M}^{-}\right) ,  \label{fordomainatrac}
\end{equation}
where $\gamma _{1},\ldots ,\gamma _{n}$ is taken from any reduced expression
$w^{-}w=s_{n}\cdots s_{1}$.
\end{proposicao}

Now define
\begin{equation}
{\mathbb S}_{\Theta }^{w}=Q_{\phi }\cdot \mathrm{st}_{\Theta }(H,w)=\{q\cdot
b\in \Bbb{E}_{\Theta }:q\in Q_{\phi },b\in \mathrm{st}_{\Theta }(H,w)\}.
\label{eq:imagem-fibrado-estavel}
\end{equation}
We plan to prove that ${\mathbb S}_{\Theta }^{w}$ is the stable set of the
Morse component $\mathcal{M}_{\Theta }^{w}$. Note that ${\mathbb S}_{\Theta
}^{w}$ is $\phi _{t}$-invariant because $Q_{\phi }$ is $\phi _{t}$-invariant
and $\mathrm{st}_{\Theta }(H,w)$ is $Z_{H}$-invariant. Also, ${\mathbb S}%
_{\Theta }^{w}$ contains $\mathcal{M}_{\Theta }^{w}$. The next result
collects some properties of these sets.

\begin{proposicao}
\label{propos:propriedades-var-estavel}Fix $\Theta \subset \Sigma $ and $%
w\in \mathcal{W}$. Then the following statements hold.

\begin{enumerate}
\item[i)]  $({\mathbb
S}_{\Theta }^{w})_{\pi (q)}=q\cdot \mathrm{st}_{\Theta }(H,w)$ for a fixed $%
q\in Q_{\phi }$.

\item[ii)]  We have that $\mathrm{cl}({\mathbb S}_{\Theta
}^{w})=\bigcup_{s\geq w}{\mathbb S}_{\Theta }^{s}$. Furthermore, if $%
\mathcal{M}_{\Theta }^{s}\subset \mathrm{cl}({\mathbb S}_{\Theta }^{w})$
then $\mathcal{M}_{\Theta }^{s}=\mathcal{M}_{\Theta }^{\overline{w}}$, where
$\overline{w}\geq w$.

\item[iii)]  The flag bundle decomposes as ${\mathbb E}_{\Theta }=\coprod \{{%
\mathbb S}_{\Theta }^{w}:\,w\in \mathcal{W}_{H}\backslash \mathcal{W}/%
\mathcal{W}_{\Theta }\}$.

\item[iv)]  In the maximal flag bundle $\mathrm{cl}({\mathbb
S}^{w})=\mathcal{A}\left( \mathcal{M}^{w}\right) =\gamma _{1}\cdots \gamma
_{n}\left( \mathcal{M}^{-}\right) $.

\item[v)]  Each ${\mathbb S}_{\Theta }^{w}$ is homeomorphic to a vector
bundle over $\mathcal{M}_{\Theta }^{w}$.
\end{enumerate}
\end{proposicao}

\begin{profe}
Fixing $q\in Q_{\phi }$ we have that $(Q_{\phi })_{\pi (q)}=qZ_{H}$ and,
since $\mathrm{st}_{\Theta }(H,w)$ is $Z_{H}$-invariant we obtain
\[
({\mathbb S}_{\Theta }^{w})_{\pi (q)}=(Q_{\phi })_{\pi (q)}\cdot \mathrm{st}%
_{\Theta }(H,w)=q\cdot \mathrm{st}_{\Theta }(H,w),
\]
which proves (i). For (ii) we take $q\in Q_{\phi }$ and prove first that
\[
\mathrm{cl}({\mathbb S}_{\Theta }^{w})_{\pi (q)}=q\cdot \mathrm{cl}(\mathrm{%
st}_{\Theta }(H,w)).
\]
The inclusion ``$\supset $'' is immediate. For the other inclusion take $\xi
\in \mathrm{cl}({\mathbb S}_{\Theta }^{w})$ such that $\pi (\xi )=\pi
(q)=x\in X$. Then there is a net $\xi _{k}\in {\mathbb S}_{\Theta }^{w}$
such that $\xi _{k}\to \xi $ so that $x_{k}:=\pi (\xi _{k})\to \pi (\xi )=x$%
. Let $\chi :U\to Q_{\phi }$ be a continuous local section of $Q_{\phi
}\subset Q$ defined in a neighborhood $U$ of $x$ and take $k$ sufficiently
large so that $x_{k}\in U$. From item (i) it follows that $\xi _{k}=\chi
(x_{k})\cdot v_{k}$ where $v_{k}\in \mathrm{st}_{\Theta }(H,w)$, so that
\[
v_{k}=\chi (x_{k})^{-1}\xi _{k}\to \chi (x)^{-1}\xi \in \mathrm{cl}(\mathrm{%
st}_{\Theta }(H,w)),
\]
which proves that $\xi \in q\cdot \mathrm{cl}(\mathrm{st}_{\Theta }(H,w))$
and establishes the other inclusion. Using the first statement and equation (%
\ref{forschubert}) for nonregular $H$, it follows that
\[
\mathrm{cl}({\mathbb S}_{\Theta }^{w})_{\pi (q)}=q\cdot \bigcup_{\overline{w}%
\geq w}\mathrm{st}_{\Theta }(H,\overline{w})=\bigcup_{\overline{w}\geq
w}q\cdot \mathrm{st}_{\Theta }(H,\overline{w})=\bigcup_{\overline{w}\geq w}({%
\mathbb S}_{\Theta }^{\overline{w}})_{\pi (q)}.
\]
Since $q\in Q_{\phi }$ is arbitrary, this proves the first statement of
(ii). The second statement of (ii) and (iii) follow by (i) and the Bruhat
decomposition of the typical fiber ${\mathbb F}_{\Theta }$.

Statement (iv) follows from Proposition \ref{propdmatalg} and the
characterizations of the Schubert cells in Subsection \ref{subsecschub}.

In order to prove (v), we first take the $K_{H}$-reduction $R_{\phi }$ of $%
Q_{\phi }$ provided by Theorem \ref{teo:tipo-parabol-flow} (iv) and the $%
L_{H}$-invariant vector bundle $\pi _{\Theta }:V_{\Theta }^{-}(H,w)\to {%
\mathrm{fix}_{\Theta }(H,w)}$ of Section \ref{seclineariz}. Since $%
K_{H}\subset L_{H}$ we can define the associated bundle of $R_{\phi }$ given
by
\[
\mathcal{V}=R_{\phi }\times _{K_{H}}V_{\Theta }^{-}(H,w)\to X.
\]
Defining the projection $\pi _{\Theta }:\mathcal{V}\to {\mathbb
E}_{\Theta }$, $q\cdot v\mapsto q\cdot \pi _{\Theta }(v)$ we have, by of
Theorem \ref{teo:tipo-parabol-flow} (iv), that $\pi _{\Theta }(\mathcal{V})=%
\mathcal{M}_{\Theta }^{w}$. This shows that $\mathcal{V}$ can be viewed as a
vector bundle over $\mathcal{M}_{\Theta }^{w}$. Recall the $K_{H}$%
-equivariant diffeomorphism $\psi :V_{\Theta }^{-}(H,w)\to \mathrm{st}%
_{\Theta }(H,w)$ of Section \ref{seclineariz} and consider the homeomorphism
\[
\Psi :\mathcal{V}\to {\mathbb E}_{\Theta },\qquad q\cdot v\mapsto q\cdot
\psi (v),
\]
whose image is precisely $R_{\phi }\cdot \mathrm{st}_{\Theta }(H,w)$. Using
that $Q_{\phi }=R_{\phi }Z_{H}$ and that $\mathrm{st}_{\Theta }(H,w)$ is $%
Z_{H}$-invariant it follows that $\Psi $ has image ${\mathbb
S}_{\Theta }^{w}$.
\end{profe}

\begin{teorema}
\label{teo:var-estavel}The stable set of the Morse component $\mathcal{M}%
_{\Theta }^{w}\subset \Bbb{E}_{\Theta }$ is
\[
\mathrm{st}(\mathcal{M}_{\Theta }^{w})={\mathbb S}_{\Theta }^{w},
\]
for every $w\in \mathcal{W}$, $\Theta \subset \Sigma $. In particular, each $%
\mathrm{st}(\mathcal{M}_{\Theta }^{w})$ is homeomorphic to a vector bundle
over $\mathcal{M}_{\Theta }^{w}$.
\end{teorema}

\begin{profe}
We work first in the maximal flag bundle ${\mathbb E}$ where we omit the
subscript $\Theta $. In this case the proof is by induction on the length of
$w$.

First we observe that if $w$ has maximal length, that is, $w=w^{-}$ then $%
\mathcal{M}^{w}$ is the repeller component by Theorem \ref
{teo:tipo-parabol-flow}. So that $\mathrm{st}(\mathcal{M}^{w})=\mathcal{M}%
^{w}={\mathbb S}^{w}$, since $\mathrm{st}_{\Theta }(H,w^{-})=w^{-}b$. Now
suppose that the result is true for all $s\in \mathcal{W}$ with length
strictly greater then the length of $w$. Then $\mathrm{st}(\mathcal{M}^{s})={%
\mathbb S}^{s}$ for all $s>w$ by the inductive hypothesis.

By (iv) and (ii) of the above proposition the domain of attraction $\mathcal{%
A}\left( \mathcal{M}^{w}\right) $ of $\mathcal{M}^{w}$ is the invariant set
\[
\bigcup_{s\geq w}{\mathbb S}^{s}.
\]
By the definitions, the stable set of $\mathcal{M}^{w}$ is contained in $%
\mathcal{A}\left( \mathcal{M}^{w}\right) $. It follows that if we take
\[
x\in \mathcal{A}\left( \mathcal{M}^{w}\right) \setminus {\mathbb
S}^{w}=\bigcup_{s>w}{\mathbb S}^{s},
\]
then its $\omega $-limit set is contained in $\mathcal{M}^{s}\neq \mathcal{M}%
^{w}$ with $s>w$. Hence the stable set of $\mathcal{M}^{w}$ is contained in $%
{\ \mathbb S}^{w}$. Since $\mathcal{A}\left( \mathcal{M}^{w}\right) $ is
closed and invariant, it contains the $\omega $-limit sets of its points, so
that the stable set of $\mathcal{M}^{w}$ is indeed ${\mathbb S}^{w}$.

Now let $\Theta \subset \Sigma $ be arbitrary and let $\pi _{\Theta }:{%
\mathbb E}\to {\mathbb E}_{\Theta }$ be the natural projection between these
flag subbundles. Then
\[
{\mathbb S}_{\Theta }^{w}=\pi _{\Theta }({\mathbb S}^{w})=\pi _{\Theta }(%
\mathrm{st}(\mathcal{M}^{w}))\subset \mathrm{st}(\mathcal{M}_{\Theta }^{w}),
\]
where the last inclusion follows from the equivariance of the flow $\phi
_{t} $ with respect to the projection $\pi _{\Theta }$. Since both ${\mathbb %
S}_{\Theta }^{w}$ and $\mathrm{st}(\mathcal{M}_{\Theta }^{w})$ partition ${%
\mathbb E}_{\Theta }$ into disjoint sets this inclusion implies the equality
${\mathbb S}_{\Theta }^{w}=\mathrm{st}(\mathcal{M}_{\Theta }^{w})$, which
proves the result.
\end{profe}

\begin{remarque}
Note that, in the present situation, the stable set of $\mathcal{M}^{w}$ is
open and dense in its domain of attraction $\mathcal{A}\left( \mathcal{M}%
^{w}\right) $. These results on stable bundles extend to flows on flag
bundles the Bruhat decomposition for the action of semi-simple elements in
flag manifolds (see e.g. \cite{dkv}).
\end{remarque}

Defining ${\mathbb U}_{\Theta }^{w}=Q_{\phi }\cdot \mathrm{un}_{\Theta
}(H,w) $, it follows by entirely analogous arguments that
\[
\mathrm{un}(\mathcal{M}_{\Theta }^{w})={\mathbb U}_{\Theta }^{w}.
\]

\section{Linearizations around Morse components\label{seclinear}}

Now we write down linearizations of the flow around the Morse components on
the flag bundles. Take a characteristic element $H_{\phi }$ for the flow and
recall the subgroup $L_{H_{\phi }}=K_{H_{\phi }}A_{H_{\phi }}\subset
Z_{H_{\phi }}$ of Section \ref{seclineariz}.

\begin{defi}
We say that the flow $\phi _{t}$ admits a conformal reduction if there
exists a subbundle $C_{\phi }\subset Q_{\phi }$ with structural group $L_{H}$
which is $\phi _{t}$-invariant.
\end{defi}

In general a flow does not admit a conformal reduction contrary to the $%
Z_{H} $-block reduction $Q_{\phi }$ that holds by Theorem \ref
{teo:tipo-parabol-flow}. To construct linearizations on the flag bundles we
make the additional assumption that conformal reduction exists. The reason
for this asumption is that linearization on flag manifolds built in
Subsection \ref{subseclinflag} is $L_{H}$-equivariant but not in general $%
Z_{H}$-equivariant. We note however that if the flow is regular, then $%
L_{H_{\phi }}=Z_{H_{\phi }}=MA$ and the conformal reduction is granted by
the block reduction. This regularity of the flow happens, for example, when
the flow has simple Lyapunov spectrum for every ergodic measure (as follows
by the results of \cite{ll1}).

When there exists a conformal reduction $C_{\phi }$ it can replace the block
reduction $Q_{\phi }$ in Theorem \ref{teo:tipo-parabol-flow}. For instance,
the Morse components are given by
\begin{equation}
{\mathcal{M}}_{\Theta }^{w}=C_{\phi }\cdot {\mathrm{fix}}_{\Theta }(H{,w)}%
=C_{\phi }\cdot wb_{\Theta }  \label{forcompredconfor}
\end{equation}
which is contained in its stable bundle
\begin{equation}  \label{eqconfst}
{\mathbb S}_{\Theta}^w = C_\phi \cdot \mathrm{st}_\Theta(H,w).
\end{equation}
In fact, ${\mathcal{M}}_{\Theta }^{w}=Q_{\phi }\cdot {\mathrm{fix}}_{\Theta
}(H,w)$ by Theorem \ref{teo:tipo-parabol-flow} (i). But then it follows from
the $Z_{H}$-invariance of ${\mathrm{fix}_{\Theta }}(H,w)$ and the equality $%
Q_{\phi }=C_{\phi }\cdot Z_{H}$ that
\[
{\mathcal{M}}_{\Theta }^{w}=C_{\phi }Z_{H}\cdot {\mathrm{fix}}_{\Theta
}(H,w)=C_{\phi }\cdot {\mathrm{fix}_{\Theta }}(H,w).
\]
Also, ${\mathrm{fix}}_{\Theta }\left( H,w\right) =L_{H}wb_{\Theta }$ so that
${\mathcal{M}}_{\Theta }^{w}=C_{\phi }\cdot wb_{\Theta }$. For the equality
involving ${\mathbb S}_{\Theta}^w$, using the equation (\ref
{eq:imagem-fibrado-estavel}), we proceed in the same way as in the proof of
the first equality above, since $\mathrm{st}_\Theta(H,w)$ is $Z_H$ invariant.

The construction of the linearization on the flag bundles follows easily
from the equivariance of the linarization on the flag manifolds. Assume
there exists a conformal reduction $C_{\phi }$ and fix $w\in \mathcal{W}$, $%
\Theta \subset \Sigma $. Recall the $L_{H}$-invariant vector bundle $\pi
_{\Theta }:V_{\Theta }(H,w)\to {\mathrm{fix}_{\Theta }(H,w)}$ of Section \ref
{seclineariz} and consider the associated bundle of $C_{\phi }$ given by
\[
\mathcal{V}_{\Theta }^{w}=C_{\phi }\times _{L_{H}}V_{\Theta }(H,w)\to X,
\]
obtained by the left action of $L_{H}$ on $V_{\Theta }(H,w)$. The bundle $%
\mathcal{V}_{\Theta }^{w}$ is a vector bundle over $\mathcal{M}_{\Theta
}^{w} $ with projection
\[
\pi _{\Theta }:\mathcal{V}_{\Theta }^{w}\to {\mathbb E}_{\Theta },\qquad
q\cdot v\mapsto q\cdot \pi _{\Theta }(v).
\]
(Note that by (\ref{forcompredconfor}) $q\cdot \pi _{\Theta }(v)\in \mathcal{%
M}_{\Theta }^{w}$ if $q\in C_{\phi }$.

Now, recall the $L_{H}$-equivariant map $\psi $ of Section \ref{seclineariz}
and consider the bundle map
\begin{equation}
\Psi :\mathcal{V}_{\Theta }^{w}\to {\mathbb E}_{\Theta },\qquad q\cdot
v\mapsto q\cdot \psi (v).  \label{eq:def-linearizacao-fibrado}
\end{equation}
It is well defined and its image contains $\mathcal{M}_{\Theta }^{w}$. Also,
$\phi _{t}$ induces a flow on associated byndle $\mathcal{V}_{\Theta }^{w}$
by
\[
\Phi _{t}(q\cdot v)=\phi _{t}(q)\cdot v.
\]
Since the following diagram commutes
\begin{equation}
\xymatrix{ \mathcal{V}_{\Theta}^w \ar[d]_{\pi_{\Theta}} \ar[rr]^{\Phi_t}& &
\mathcal{V}_{\Theta}^w \ar[d]^{\pi_{\Theta}} \\ \mathcal{M}_{\Theta}^w
\ar[rr]_{\phi_t} & & \mathcal{M}_{\Theta}^w}
\end{equation}
it follows that $\Phi _{t}$ is a linear flow induced by $\phi _{t}$ such
that its induced flow on the base $\mathcal{M}_{\Theta }^{w}$ is precisely
the restriction of $\phi _{t}$, which is chain transitive. By of Theorem \ref
{teolineariz} (i) it follows that $\Phi _{t}$ and $\phi _{t}$ are conjugated
by the equivariant map $\Psi $ in a neighborhood of the null section.

We define the associated bundle of $C_\phi$ given by
\begin{equation}  \label{eqstable}
\mathcal{S}_{\Theta}^w = C_\phi \times_{L_H} V^-_\Theta(H,w),
\end{equation}
which, in the same way as above, can be viewed as a vector bundle over $%
\mathcal{M}^w_\Theta$. By the equation (\ref{eqconfst}), its image by $\Psi$
is precisely the $\phi_t$-invariant set ${\mathbb
S}_{\Theta}^w$ which contains $\mathcal{M}^w_\Theta$. By item (ii) of the
Theorem \ref{teolineariz}, it follows that $\Phi_t$ and $\phi_t$ are
conjugated by the equivariant global homeomorphism $\Psi$ restricted to the $%
\phi_t$-invariant vector bundle $\mathcal{S}_{\Theta}^w$, which from now on
will be called the stable vector bundle of $\mathcal{M}_{\Theta}^w$. We
collect the above results in the following proposition.

\begin{proposicao}
\label{prop3lineariz}Assume the existence of a conformal reduction $C_{\phi
} $. Then the restriction of $\phi _{t}$ to ${\mathbb S}_{\Theta }^{w}$ is
conjugate to a linear flow over the restriction of $\phi _{t}$ to $\mathcal{M%
}_{\Theta }^{w}$.
\end{proposicao}

In order to define the unstable vector bundle we proceed analogously.
Recalling the $L_H$-invariant subbundle $V^+_\Theta(H,w)$ of Section \ref
{seclineariz}, we define the $\phi_t$-invariant bundle of $C_\phi$ given by
\[
\mathcal{U}_{\Theta}^w = C_\phi \times_{L_H} V^+_\Theta(H,w).
\]
Its image under $\Psi$ is $\mathrm{un}(\mathcal{M}^w_\Theta) = {\mathbb U}%
_{\Theta}^w$. We observe that from the definition of these vector bundles it
follows from equation (\ref{eqwhitneysum}) the following $\phi_t$-invariant
Whitney sum decomposition
\begin{equation}  \label{eqdynwhitneysum}
\mathcal{V}_{\Theta}^w = \mathcal{S}_{\Theta}^w \oplus \mathcal{U}%
_{\Theta}^w \to \mathcal{M}^w_\Theta.
\end{equation}

A case where the assumption of a conformal reduction is not needed is the
following. If $\Theta \subset \Theta (\phi )$ or $\Theta (\phi )\subset
\Theta $ then, by Corollary \ref{prop1lineariz}, the conjugation $\psi $
around the attractor component is $Z_{H}$ -equivariant. Hence, proceeding as
above, we have the following.

\begin{corolario}
\label{cor3lineariz}Assume that $\Theta \subset \Theta (\phi )$ or $\Theta
(\phi )\subset \Theta $. Then the restriction of $\phi _{t}$ to ${\mathbb
S}_{\Theta }^{+}$ is conjugate to a linear flow over the restriction of $%
\phi _{t}$ to the attractor $\mathcal{M}_{\Theta }^{+}$.
\end{corolario}

\section{Conley indices\label{secconley}}

In order to compute the Conley indices associated to the finest Morse
decomposition of ${\mathbb E}_{\Theta }$, we assume, as in the previous
section, the existence of a conformal reduction $C_{\phi }$ . This condition
will then be dropped for the computation of the Conley index of the
attractor component.

\begin{proposicao}
\label{propnorm} Assume the existence of a conformal reduction $C_{\phi }$.
Then there exists a norm $|\cdot |$ in the vector bundle $\mathcal{V}%
_{\Theta }^{w}$ and a constant $\alpha >0$ such that

\begin{itemize}
\item[i)]  $|\Phi _{t}(q\cdot v)|\leq e^{-\alpha t}|q\cdot v|$, for $q\cdot
v\in \mathcal{S}_{\Theta }^{w}$,

\item[ii)]  $|\Phi _{t}(q\cdot v)|\geq e^{\alpha t}|q\cdot v|$, for $q\cdot
v\in \mathcal{U}_{\Theta }^{w}$.
\end{itemize}
\end{proposicao}

\begin{profe}
Let $(\mathcal{S}_{\Theta }^{w})_{0}$ be the null section of the stable
vector bundle $\mathcal{S}_{\Theta }^{w}\to \mathcal{M}_{\Theta }^{w}$
defined in the equation (\ref{eqstable}). Using the linearization $\Psi $
given by Proposition \ref{prop3lineariz} and using Theorem \ref
{teo:var-estavel}, it follows that, for each $v\in \mathcal{S}_{\Theta }^{w}$%
, we have $\Phi _{t}(v)\to (\mathcal{S}_{\Theta }^{w})_{0}$ when $t\to
\infty $. Choosing a continuous vector bundle norm $|\cdot |_{1}$ in $%
\mathcal{S}_{\Theta }^{w}\to \mathcal{M}^{w}$ we get from Fenichel's
uniformity Lemma (see e.g. \cite{ayala2}, Proposition 3.2) that there exists
constants $C,\alpha _{1}>0$ such that
\[
|\Phi _{t}(q\cdot v)|_{1}\leq Ce^{-\alpha _{1}t}|q\cdot v|_{1},
\]
for $q\cdot v\in \mathcal{S}_{\Theta }^{w}$. By Proposition 3.2 of \cite
{ayala2} we can change the above norm to obtain $C=1$ in the above
inequality, so we can assume without lost of generality that $C=1$. In an
entirely analogous fashion, arguing with the reversed time flow in the
unstable vector bundle $\mathcal{U}_{\Theta }^{w}\to \mathcal{M}_{\Theta
}^{w}$, we get a norm $|\cdot |_{2}$ in $\mathcal{U}_{\Theta }^{w}$ and a
constant $\alpha _{2}>0$ such that
\[
|\Phi _{t}(q\cdot v)|_{2}\geq e^{\alpha _{2}t}|q\cdot v|_{2},
\]
for $q\cdot v\in \mathcal{U}_{\Theta }^{w}$. Using the Whitney sum
decomposition $\mathcal{V}_{\Theta }^{w}=\mathcal{S}_{\Theta }^{w}\oplus
\mathcal{U}_{\Theta }^{w}\to \mathcal{M}_{\Theta }^{w}$ given in (\ref
{eqdynwhitneysum}) it is clear that these two norms add up to a norm in $%
|\cdot |$ in $\mathcal{V}_{\Theta }^{w}$ which restricts to $|\cdot |_{1}$
in $\mathcal{S}_{\Theta }^{w}$ and to $|\cdot |_{2}$ in $\mathcal{U}_{\Theta
}^{w}$. Choosing $\alpha :=\min \{\alpha _{1},\,\alpha _{2}\}$ we obtain the
result.
\end{profe}

The previous result shows that the linearized flow in $\mathcal{V}%
_{\Theta}^w $ is hyperbolic, so we can think of each Morse component $%
\mathcal{M}^w_\Theta$ as a normally hyperbolic invariant set in ${\mathbb E}%
_\Theta$. For the attractor component, using Corollary \ref{cor3lineariz}
and proceeding as in the above proof, we can drop the assumption on the
existence of a conformal reduction.

\begin{corolario}
\label{coratrator} Assume that $\Theta \subset \Theta (\phi )$ or $\Theta
(\phi )\subset \Theta $. Then there exists a norm $|\cdot |$ in the vector
bundle $\mathcal{S}_{\Theta }^{+}$ and a constant $\alpha >0$ such that $%
|\Phi _{t}(q\cdot v)|\leq e^{-\alpha t}|q\cdot v|$, for $q\cdot v\in
\mathcal{S}_{\Theta }^{+}$.
\end{corolario}

Now we use the linearization $\Psi :\mathcal{V}_{\Theta }^{w}\to {\mathbb E}%
_{\Theta }$ given in (\ref{eq:def-linearizacao-fibrado}) and the norm in $%
\mathcal{V}_{\Theta }^{w}$ given by Proposition \ref{propnorm} to construct
an index pair for each $\mathcal{M}_{\Theta }^{w}$. Recall that, by item (i)
of Theorem \ref{teolineariz}, it follows that $\Phi _{t}$ and $\phi _{t}$
are conjugated by the equivariant map $\Psi $ in a neighborhood $A$ of the
null section. For $q\in C_{\phi }$, $v\in V_{\Theta }(H,w)$, let $q\cdot
v^{s}$ and $q\cdot v^{u}$ denote, respectively, the stable and the unstable
components of the invariant Whitney sum decomposition of $\mathcal{V}%
_{\Theta }^{w}$ given in (\ref{eqdynwhitneysum}). We can multiply the norm
of $\mathcal{V}_{\Theta }^{w}$ by a positive constant so that the set $%
\{q\cdot v:|q\cdot v^{s}|,|q\cdot v^{u}|\leq 1\}$ is contained in $A$. We
then define
\[
{\mathbb B}_{1}=\Psi (\{q\cdot v:|q\cdot v^{s}|,|q\cdot v^{u}|\leq 1\})
\]
and
\[
{\mathbb B}_{0}=\Psi (\{q\cdot v:|q\cdot v^{s}|\leq 1,|q\cdot v^{u}|=1\}).
\]

\begin{proposicao}
Assume the existence of a conformal reduction $C_{\phi }$. Then the pair $({%
\mathbb B}_{1},{\mathbb B}_{0})$ is an index pair for $\mathcal{M}_{\Theta
}^{w}$.
\end{proposicao}

\begin{profe}
We need to verify the three properties of the definition of index pair (see
Section \ref{preliminar-conley}). Property (1) is immediate, since clearly $%
\mathrm{cl}({\mathbb B}_{1}\backslash {\mathbb B}_{0})={\mathbb B}_{1}$ is a
neighborhood of $\mathcal{M}_{\Theta }^{w}$ which, from Proposition \ref
{propnorm}, contains no other invariant sets. In order to check property
(2), we first note that
\[
\phi _{t}(\Psi (q\cdot v))=\Psi (\Phi _{t}(q\cdot v)).
\]
Then, for $\xi \in {\mathbb B}_{0}$ and $t\geq 0$, we have $\phi _{t}(\xi
)\in {\mathbb B}_{1}$ if and only if $t=0$, since
\[
|\Phi _{t}(q\cdot v^{u})|\geq e^{\alpha t}|q\cdot v^{u}|=e^{\alpha t}\geq 1,
\]
by Proposition \ref{propnorm}. It remains to prove property (3). Thus let $%
\xi \in {\mathbb B}_{1}$ be such that there exists $t^{\prime }>0$ with $%
\phi _{t^{\prime }}(\xi )\notin {\mathbb B}_{1}$. We define
\[
t_{m}=\sup \{t:\phi _{t}(\xi )\in {\mathbb B}_{1}\}.
\]
Since ${\mathbb B}_{1}$ is closed, it follows that $\phi _{t_{m}}(\xi )\in {%
\mathbb B}_{1}$. We must show that $\phi _{t_{m}}(\xi )\in {\mathbb B}_{0}$.
If we assume the contrary then $\phi _{t_{m}}(\xi )=\Psi (q\cdot v)$ with $%
|q\cdot v^{s}|\leq 1$ and $|q\cdot v^{u}|<1$. For $t>0$ we have
\[
\phi _{t_{m}+t}(\xi )=\phi _{t}(\Psi (q\cdot v))=\Psi (\Phi _{t}(q\cdot v)).
\]
But
\[
|\Phi _{t}(q\cdot v^{s})|\leq e^{-\alpha t}|q\cdot v^{s}|<1.
\]
Hence by continuity we also have $|\Phi _{t}(q\cdot v^{u})|<1$ if $t>0$
small enough. This implies that $\phi _{t_{m}+t}(\xi )\in {\mathbb B}_{1}$,
contradicting the definition of $t_{m}$.
\end{profe}

Now we relate the Conley index of $\mathcal{M}_{\Theta }^{w}$ with the Thom
space (see Section \ref{preliminar-conley}) of the unstable vector bundle $%
\mathcal{U}_{\Theta }^{w}\to \mathcal{M}_{\Theta }^{w}$. This relation is
well know for hyperbolic linear flows.

\begin{teorema}
\label{teo:indice-conley}Assume the existence of a conformal reduction $%
C_{\phi }$. Then the Conley index of $\mathcal{M}_{\Theta }^{w}$ is the
homotopy class of the Thom space of the unstable vector bundle $\mathcal{U}%
_{\Theta }^{w}\to \mathcal{M}_{\Theta }^{w}$.
\end{teorema}

\begin{profe}
We start the proof by defining
\[
{\mathbb B}_{1}^{u}=\{\Psi (q\cdot v^{u}):|q\cdot v^{u}|\leq 1\}\subset {%
\mathbb B}_{1}
\]
and
\[
{\mathbb B}_{0}^{u}=\{\Psi (q\cdot v^{u}):|q\cdot v^{u}|=1\}\subset {\mathbb %
B}_{0}.
\]
We note that the inclusion $i:{\mathbb B}_{1}^{u}\to {\mathbb B}_{1}$ and
the deformation retract
\[
r:{\mathbb B}_{1}\to {\mathbb B}_{1}^{u},\qquad q\cdot v\mapsto q\cdot v^{u},
\]
are such that $r\circ i=\mathrm{id}_{{\mathbb B}_{1}^{u}}$ and, clearly,
there exists a continuous homotopy $h:[0,1]\times {\mathbb B}_{1}\to {%
\mathbb B}_{1}$ between the map $i\circ r$ and the map $\mathrm{id}_{{%
\mathbb B}_{1}}$ such that
\[
h([0,1]\times {\mathbb B}_{0})\subset {\mathbb B}_{0}.
\]

Since the Conley index is given by the homotopy class of the quotient space $%
{\mathbb B}_{1}/{\mathbb B}_{0}$ and the Thom space of the unstable vector
bundle ${\mathbb U}_{\Theta }^{w}$ is the quotient space ${\mathbb B}%
_{1}^{u}/{\mathbb B}_{0}^{u}$, we need to show that these quotient spaces
are homotopically equivalent. Since $i({\mathbb B}_{0}^{u})\subset {\mathbb B%
}_{0}$ and $r({\mathbb
B}_{0})={\mathbb B}_{0}^{u}$, we can define maps $I$ and $R$ such that the
following diagrams commute
\begin{equation}
\xymatrix{ {\mathbb B}_1 \ar[d]_{p} \ar[rr]^{r}& & {\mathbb B}_1^u
\ar[d]_{p^u} \ar[rr]^{i}& & {\mathbb B}_1 \ar[d]^{p} \\ {\mathbb B}_1 /
{\mathbb B}_0 \ar@{-->}[rr]_{R} & & {\mathbb B}_1^u / {\mathbb B}_0^u
\ar@{-->}[rr]_{I} & & {\mathbb B}_1 / {\mathbb B}_0}
\end{equation}
where $p$ and $p^{u}$ are the respective projections. These diagrams imply
that $I$ and $R$ are continuous maps and it is immediate to show that $%
R\circ I=\mathrm{id}_{{\mathbb B}_{1}^{u}/{\mathbb
B}_{0}^{u}}$. Finally, since
\[
h([0,1]\times {\mathbb B}_{0})\subset {\mathbb B}_{0},
\]
we can define a map $H$ such that the following diagram commutes
\begin{equation}
\xymatrix{ [0,1] \times {\mathbb B}_1 \ar[d]_{{\rm{id}}_{[0,1]} \times p}
\ar[rr]^{h}& & {\mathbb B}_1 \ar[d]^{p} \\ [0,1] \times ({\mathbb B}_1 /
{\mathbb B}_0) \ar@{-->}[rr]_{H} & & {\mathbb B}_1 / {\mathbb B}_0}
\end{equation}
showing that $H$ is a continuous homotopy between the map $I\circ R$ and the
map $\mathrm{id}_{{\mathbb B}_{1}/{\mathbb B}_{0}}$ and therefore that ${%
\mathbb B}_{1}/{\mathbb B}_{0}$ and ${\mathbb
B}_{1}^{u}/{\mathbb B}_{0}^{u}$ are homotopically equivalent spaces.
\end{profe}

\begin{corolario}
The following isomorphism in cohomology holds
\begin{equation}
CH^{*+n_{w}}(\mathcal{M}_{\Theta }^{w})\simeq H^{*}(\mathcal{M}_{\Theta
}^{w}),  \label{isom:indice-conley}
\end{equation}
where $n_{w}$ is the dimension of $\mathcal{U}_{\Theta }^{w}$. The
cohomology coefficients are taken in ${\mathbb Z}_{2}$ in the general case
and in ${\mathbb Z}$ if $\mathcal{U}_{\Theta }^{w}$ is orientable.
\end{corolario}

\begin{profe}
Follows directly from Thom's isomorphism (see Section \ref{preliminar-conley}%
).
\end{profe}

Since there is also a Thom isomorphism in singular homology, the above
isomorphisms also hold for the singular homology of the Conley index.

We now specialize the previous result to the situation of a contractible
base (which is related to control flows).

\begin{proposicao}
\label{prop:indice-conley-contratil}Assume the existence of a conformal
reduction $C_{\phi }$. Put $\Delta =\Theta (H)$ and take $H_{\Theta }\in
\mathrm{cl}\frak{a}^{+}$ such that $\Theta =\Theta (H_{\Theta })$. If the
base $X$ is contractible, then the Conley index of $\mathcal{M}_{\Theta }^{w}
$ is the homotopy class of the Thom space of the vector bundle $V_{\Theta
}^{+}(H,w)\to \mathbb{F}(\Delta)_{H_{0}}$, where $H_{0}$ is the orthogonal
projection of $wH_{\Theta }$ in $\frak{a}\left( \Delta \right) $. In
particular, we have the following isomorphism in cohomology
\[
CH^{*+n_{w}}(\mathcal{M}_{\Theta }^{w})\simeq H^{*}\left( \mathbb{F}(\Delta)%
_{H_{0}}\right) ,
\]
where $n_{w}$ is the dimension of $V_{\Theta }^{+}(H,w)$. The cohomology
coefficients are taken in ${\mathbb Z}_{2}$ in the general case and in ${%
\mathbb Z}$ if $V_{\Theta }^{+}(H,w)$ is orientable.
\end{proposicao}

\begin{profe}
Since $X$ is contractible  $C_{\phi }$ is trivial (see \cite{steenrod},
Corollary 11.6) thus there exists a continuos global section $\chi :X\to
C_{\phi }$. It follows that the application $\lambda :X\times V_{\Theta
}^{+}(H,w)\to \mathcal{U}_{\Theta }^{w}$, given by
\[
(x,v)\mapsto \frac{\chi (x)\cdot v}{|\chi (x)\cdot v|}\Vert v\Vert
\]
is a homeomorphism whose inverse is given by
\[
\chi (x)\cdot v\mapsto \left( x,\frac{v}{\Vert v\Vert }|\chi (x)\cdot
v|\right) ,
\]
where $|\cdot |$ is the norm in $\mathcal{U}_{\Theta }^{w}$ given by
Proposition \ref{propnorm} and $\Vert \cdot \Vert $ is some norm in $%
V_{\Theta }^{+}(H,w)$. Let
\[
B_{1}^{u}=\{v\in V_{\Theta }^{+}(H,w):\Vert v\Vert \leq 1\}
\]
and
\[
B_{0}^{u}=\{v\in V_{\Theta }^{+}(H,w):\Vert v\Vert =1\}.
\]
It is then straightforward that $\lambda (X\times B_{1}^{u})$ is contained
in the open set where $\Psi $ is an homeomorphism and also that
\[
\Psi (\lambda (X\times B_{1}^{u}))={\mathbb B}_{1}^{u}\qquad \mbox{and}%
\qquad \Psi (\lambda (X\times B_{0}^{u}))={\mathbb B}_{0}^{u}.
\]
Thus we have
\[
{\mathbb B}_{1}^{u}/{\mathbb B}_{0}^{u}\simeq (X\times B_{1}^{u})/(X\times
B_{0}^{u})
\]
whose homotopy class is, by Theorem \ref{teo:indice-conley}, the Conley
index of $\mathcal{M}_{\Theta }^{w}$. It remains to show that the space $%
(X\times B_{1}^{u})/(X\times B_{0}^{u})$ and the space $B_{1}^{u}/B_{0}^{u}$
are homotopically equivalent. We apply the  same argument of  Theorem \ref
{teo:indice-conley} where the inclusion is given by
\[
i:B_{1}^{u}\to X\times B_{1}^{u},\qquad v\mapsto (\overline{x},v),
\]
the deformation retract is given by
\[
r:X\times B_{1}^{u}\to B_{1}^{u},\qquad (x,v)\mapsto v
\]
and the homotopy is given by
\[
h:[0,1]\times (X\times B_{1}^{u})\to X\times B_{1}^{u},\qquad
(t,(x,v))\mapsto (g(t,x),v)
\]
where $g:[0,1]\times X\to X$ is some homotopy between $\mathrm{id}_{X}$ and
the constant map $\overline{x}\in X$.
\end{profe}

Now we apply  the abstract Morse equation to our flows. By Theorem \ref
{teo:tipo-parabol-flow}, the finest Morse decomposition in ${\mathbb E}%
_{\Theta }$ is given by $\mathcal{M}_{\Theta }^{w}$, for $w\in \mathcal{W}%
_{H}\backslash \mathcal{W}/\mathcal{W}_{\Theta }$. By Theorem \ref
{teo:indice-conley} %, which applies to $\check{{\rm C}}$ech cohomology,
we have that $CH^{*+n_{w}}(\mathcal{M}_{\Theta }^{w})\simeq H^{*}(\mathcal{M}%
_{\Theta }^{w})$, where $n_{w}$ is the dimension of the unstable bundle of $%
\mathcal{M}^{w}$. In the maximal flag bundle ${\mathbb E}$, item (ii) of
Proposition \ref{prop:compon-subflags} gives us that $H^{*}(\mathcal{M}%
^{w})\simeq H^{*}(\mathcal{M}^{+})$. It follows that
\[
CP(t,\mathcal{M}^{w})=\sum_{j\geq 0}t^{j+n_{w}}H^{j}(\mathcal{M}%
^{w})=t^{n_{w}}P(t,\mathcal{M}^{w})=t^{n_{w}}P(t,\mathcal{M}^{+})
\]
By Theorem \ref{teo:morse-eq-abstrata} we have then immediately.

\begin{corolario}
Assume the existence of a conformal reduction $C_{\phi }$ and suppose that
the cohomology of the base $X$ have finite rank in all dimensions. Then the
Morse equation in the maximal flag bundle ${\mathbb E}$ becomes
\begin{equation}
P(t,\mathcal{M}^{+})\sum_{w\in \mathcal{W}_{H}\backslash \mathcal{W}%
}t^{n_{w}}=P(t,{\mathbb E})+(1+t)R(t).  \label{formorseequation}
\end{equation}
If the coefficient ring is the integers $\Bbb{Z}$ then the coefficients of $%
R(t)$ are non-negative.
\end{corolario}

We note that where $\mathcal{M}^{+}\to X$ is a bundle with fiber ${\mathbb
F}(\Theta (\phi ))$, where $\Theta (\phi )$ is the parabolic type of the
flow $\phi _{t}$.

Equation (\ref{formorseequation}) relates the topology of the attractor
component and the unstable dimensions with the topology of the whole maximal
flag bundle ${\mathbb E}$. In the other flag bundles ${\mathbb
E}_{\Theta }$ the Morse components $\mathcal{M}_{\Theta }^{w}$ are not
necessarily homeomorphic, so that the Morse equation is not simple  as
above. Nevertheless, if the flow is regular, then each Morse component $%
\mathcal{M}_{\Theta }^{w}$ in ${\mathbb E}_{\Theta }$ is homeomorphic to the
base space $X$ so that the Morse equation for the flow in ${\mathbb
E}_{\Theta }$ becomes
\[
P(t,X)=\sum \{t^{n_{w}}:\,w\in {\mathcal{W}}/{\mathcal{W}}_{\Theta }\}=P(t,{%
\mathbb
E}_{\Theta })+(1+t)R(t).
\]
Note that in the maximal flag bundle we have $\mathcal{W}_{\Theta }=\{1\}$.

For the attractor component we can drop  the assumption on the existence of
a conformal reduction.

\begin{corolario}
The cohomological Conley index of the attractor in ${\mathbb
E}_{\Theta }$ is  given by
\[
CH^{*}(\mathcal{M}_{\Theta }^{+})=H^{*}(\mathcal{M}_{\Theta }^{+}),
\]
for every $\Theta \subset \Sigma $. Furthermore, if $\Theta \subset \Theta
(\phi )$ or $\Theta (\phi )\subset \Theta $, then the Conley index of $%
\mathcal{M}_{\Theta }^{+}$ is the homotopy class of the flag bundle $%
\mathcal{M}_{\Theta }^{+}\to X$.

Put $\Delta =\Theta (H)$ and take $H_{\Theta }\in \mathrm{cl}\frak{a}^{+}$
such that $\Theta =\Theta (H_{\Theta })$. If the base $X$ is contractible
then  the cohomological Conley index of the attractor in ${\mathbb E}%
_{\Theta }$ is  given by
\[
CH^{*}(\mathcal{M}_{\Theta }^{+})=H^{*}(\mathbb{F}(\Delta)_{H_{0}}),
\]
where $H_{0}$ is the orthogonal projection of $H_{\Theta }$ in $\frak{a}%
\left( \Delta \right) $. Also, if $\Theta \subset \Theta (\phi )$ or $\Theta
(\phi )\subset \Theta $, then the Conley index of $\mathcal{M}_{\Theta }^{+}$
is the homotopy class of the flag manifold $\mathbb{F}(\Delta)_{H_{0}}$.
\end{corolario}

\begin{profe}
The first assertion follows directly from Corollary 5 of \cite{sanj}. The
second assertion follows from Corollary \ref{coratrator}, proceeding as in
the proof of Theorem \ref{teo:indice-conley}. The third and fourth
assertions follow from the first and second assertions, item (i) of
Proposition \ref{prop:compon-subflags} and Corollary 11.6 of \cite{steenrod}%
, since $X$ is contractible.
\end{profe}

\end{document}